\documentclass[12pt]{article}
\usepackage[utf8]{inputenc}
\usepackage[T1]{fontenc}
\usepackage[cal=boondox,calscaled=.96]{mathalfa}
\usepackage{euscript}
\usepackage[debug]{hyperref}
\usepackage{tikz-cd}
\usepackage{amscd,amsfonts,amsmath,amssymb,amstext,amsthm,dsfont,mathtools}
\usepackage{latexsym}
\usepackage{enumerate}
\usepackage[all]{xy}
\usepackage{mathabx}

\newcommand{\CC}{\mathbb{C}}
\newcommand{\RR}{\mathbb{R}}
\newcommand{\NN}{\mathbb{N}}

\newcommand{\ZZ}{\mathbb{Z}}

\newcommand{\Lotimes}{\stackrel{\mathbb{L}}{\otimes}}

\newcommand{\ubi}{{\underline{{\bf i}}}}
\newcommand{\ubj}{{\underline{{\bf j}}}}

\DeclareMathOperator{\Hom}{\rm Hom}

\DeclareMathOperator{\wt}{\rm wt}

\DeclareMathOperator{\Ann}{\rm Ann}

\DeclareMathOperator{\ann}{\rm ann}

\DeclareMathOperator{\gr}{\rm gr}

\DeclareMathOperator{\Der}{\rm Der}

\DeclareMathOperator{\fDer}{\mathcal{D}\mathit{er}}

\DeclareMathOperator{\fHom}{\mathcal{H}\!\mathit{om}}
\DeclareMathOperator{\fExt}{\mathcal{E}\!\mathit{xt}}

\DeclareMathOperator{\Sym}{{\rm Sym}}

\DeclareMathOperator{\DR}{DR}

\newcommand{\OO}{\mathcal{O}}
\newcommand{\DD}{\mathcal{D}}
\newcommand{\VV}{\mathcal{V}}
\newcommand\calV{\mathcal{V}}

\newcommand{\MM}{\mathcal{M}}

\newcommand{\cM}{\mathcal{M}}

\newcommand{\calE}{\mathcal{E}}

\newcommand{\cI}{\mathcal{I}}
\newcommand{\cJ}{\mathcal{J}}

\DeclareMathOperator{\SP}{\rm Sp}
\newcommand{\dsD}{\mathds{D}}

\newcommand{\Spb}{{\mathrm{Sp}}^\bullet} 

\newcommand{\Sp}{{\mathrm{Sp}}}
\newcommand{\ord}{{\mathrm{ord}}}

\newcommand{\cD}{\DD}
\newcommand{\cO}{\OO}
\newcommand\Char{\text{\rm Char}}

\newcommand\bu{\bullet}

\newcommand\p{{\partial}}

\newcommand\vsn{\vskip 10pt\noindent}

\newcommand\ssm{\smallsetminus}

\newcommand\beq{\begin{equation}}
\newcommand\eeq{\end{equation}}
\newcommand\bc{\begin{center}}
\newcommand\ec{\end{center}}
\newcommand\ld{{\ldots}}

\newcommand\depth{{\mbox{depth}}}

\newcommand\G{{\Gamma}}

\DeclareMathOperator{\Mod}{\rm Mod}

\DeclareMathOperator{\supp}{\rm supp}

\newcounter{numero}[section]
\renewcommand{\thenumero}{(\arabic{section}.\arabic{numero})}
\newenvironment{corolario}{\medskip
\refstepcounter{numero}\noindent {\sc  \thenumero\ Corollary.}\
\it}{\vspace{1ex}\par}

\newenvironment{teorema}{\medskip
\refstepcounter{numero}\noindent {\sc  \thenumero\ Theorem.}\
\it}{\vspace{1ex}\par}

\newenvironment{lema}{\medskip
\refstepcounter{numero}\noindent {\sc  \thenumero\ Lemma.}\
\it}{\vspace{1ex}\par}

\newenvironment{definicion}{\medskip
\refstepcounter{numero}\noindent {\sc  \thenumero\ Definition.}\
\it}{\vspace{1ex}\par}

\newenvironment{proposicion}{\medskip
\refstepcounter{numero}\noindent {\sc  \thenumero\
Proposition.}\ \it}{\vspace{1ex}\par}

\newenvironment{nota}{\medskip
\refstepcounter{numero}\noindent {\sc  \thenumero\ Remark.}\
}{\vspace{1ex}\par}

\newenvironment{conjetura}{\medskip
\refstepcounter{numero}\noindent {\sc  \thenumero\ Conjecture.}\
}{\vspace{1ex}\par}

\newenvironment{ejemplo}{\medskip
\refstepcounter{numero}\noindent {\sc  \thenumero\ Example.}\
}{\vspace{1ex}\par}

\newenvironment{prueba}{
\noindent {\sc  Proof.}\ }{\hfill Q.E.D.\vspace{1ex}\par}

\newfont{\Bb}{msbm10}
\newfont{\Bc}{msbm8}

\newcommand\tos{{\ \rightarrow\ }}

\begin{document}

\title{Logarithmic Comparison Theorems}

\author{F.~J.~Castro-Jiménez\thanks{Departamento de \'Algebra \&\ Instituto de Matem\'aticas (IMUS),
Facultad de Matem\'aticas, Universidad de Sevilla, Spain,
 e-mail castro@us.es}, D.~Mond\thanks{Mathematics Institute, University of Warwick, Coventry CV4 7AL United Kingdom,
e-mail: d.m.q.mond@warwick.ac.uk}, L.~Narváez-Macarro\thanks{Departamento de \'Algebra \&\ Instituto de Matem\'aticas (IMUS),
Facultad de Matem\'aticas, Universidad de Sevilla, Spain,
 e-mail narvaez@us.es}}

\maketitle
\noindent We dedicate this chapter to the memory of our colleagues Jim Damon and José Luis Vicente Córdoba.

\begin{abstract}
In this paper we study the comparison between the logarithmic and the meromorphic de Rham complexes along a divisor in a complex manifold. We focus on the case of free divisors, starting with the case of locally quasihomogeneous divisors, and we explain how D-module theory can be used for this comparison.
\end{abstract}

\section*{Introduction}

We survey a number of results known as {\it comparison theorems}, along the lines of the comparison theorem of Grothendieck. Grothendieck's Comparison Theorem states that, for any hypersurface $D$ in a complex manifold $X$, the cohomology of the meromorphic de Rham  complex with respect to $D\subset X$ (i.e. $h^q(\Omega_X^\bullet(\star D))$ for $0 \leq q \leq \dim X$)  coincides with $\RR^q j_* \CC_{U}$ where $j : U:= X\ssm D \rightarrow X$ is the inclusion and $\CC_U$ is the constant sheaf $\CC$ on $U$. In effect, the hypercohomology of this complex is the (topological) cohomology $H^\bullet(U;\CC)$.
If $D$ is a normal crossing divisor (NCD), this comparison result was proved by Atiyah and Hodge, and Grothendieck's proof reduces the general case to the case of a NCD by using Hironaka's resolution of singularities. When $D$ is an arrangement of hyperplanes in $\CC^n$, the  Brieskorn complex $B^\bullet$ computes the singular cohomology of the complement $U=\CC^n\ssm D$, $H^\bullet(U;\CC)$; and this is also a comparison result. There is a class of divisors, the so-called {\it free divisors}, introduced by Kyoji Saito, for which a certain number of comparison results can be proven. We survey these results and show how $\cD$--module theory provides a way to deal with them.

The content of the paper is as follows. In Sect. \ref{Grothendieck-comparison-theorem} we first recall Grothendieck's comparison theorem (see also Sect. \ref{Gro-revisited}), that there is a canonical isomorphism
$$h^q(\Omega_X^\bullet(\star D)) \rightarrow \RR^q j_* \CC_{U}$$ for $0\leq q \leq \dim X$. If $X$ is Stein, then the global morphisms analogous to the previous ones $$h^q(\Gamma(X,\Omega_X^\bullet(\star D))) \rightarrow H^q(U; \CC)$$  are isomorphisms. We then recall Brieskorn's Theorem proving that if $D$ is a finite union of hyperplanes in $X=\CC^n$, with equations $h_j=0$, then the cohomology of the complement $H^q(X\ssm D;\CC)$ can be computed as the cohomology of the so called Brieskorn complex $B^\bullet$, which is the $\CC$-subalgebra of the exterior algebra $\Gamma(X,\Omega_X^\bullet(\star D))$ generated by the forms $\frac{dh_j}{h_j}$.

Kyoji Saito (see Sect. \ref{LCTfirst}) introduced the notion of {\em logarithmic meromorphic form} with respect to a divisor $D$ in a complex manifold $X$: a meromorphic form $\omega$ has a {\it logarithmic pole} along the divisor $D$ if both $h\omega$ and $hd\omega$ are holomorphic, where $h$ is a local equation for $D$. Note in particular that all of Brieskorn's forms have logarithmic poles along an arrangement of hyperplanes $D\subset \CC^n$. Logarithmic meromorphic forms form an $\cO_X$-subcomplex of the meromorphic de Rham complex $\Omega_X^\bullet(\star D)$. A natural question asks for the class of divisors $D\subset X$ for which the inclusion $$ \Omega_X^\bullet(\log D) \hookrightarrow \Omega_X^\bullet(\star D)$$ is a quasi-isomorphism (i.e. when the two complexes have the same cohomology). By analogy with Grothendieck's comparison theorem, if the morphism $ \Omega_X^\bullet(\log D) \hookrightarrow \Omega_X^\bullet(\star D)$  is a quasi-isomorphism we say that the divisor $D$ satisfies the {\em logarithmic comparison theorem} (LCT), or that LCT holds for $D$. Theorem \ref{lct} states that any {\em locally quasihomogeneous free divisor} $D \subset X$ satisfies LCT. In 1977 H. Terao conjectured (see Sect. \ref{conjecture-of-terao}) that LCT always holds for hyperplane arrangements (which are, of course, locally quasihomogeneous), and this conjecture was finally proved by Daniel Bath in \cite{Bath2022}.

A free divisor is {\em linear} if the module $\fDer(-\log D)$ of logarithmic derivations with respect to $D$ (as defined by K. Saito in \cite{Saito-80}, see Sect. \ref{LCTfirst}) has a basis of vector fields whose coefficients are linear forms. In Sect. \ref{linear-free-divisors} we recall a number of results on linear free divisors and in particular the theorem of \cite{gmns}, that LCT holds for ``reductive'' linear free divisors.

In Sect. \ref{D_X-modules} we introduce the background in $\cD$--module theory that is used in Sect. \ref{sec:free-divisors-log-D-modules}. This last section is devoted to state a $\cD$--module criterion for the Logarithmic Comparison Theorem for free divisors. First, and following the paper of F.J. Calderón-Moreno \cite{calde_ens}, we consider the sheaf $ \VV_X^D$ of logarithmic differential operators on $X$, with respect to a free divisor $D\subset X$ and the logarithmic Spencer complex $\Sp^\bullet(\log D)$ which is a locally free resolution of the $\VV_X^D$--module $\cO_X$. By using \cite{calde_nar_fourier} one generalizes the construction of this last complex to an arbitrary left $\VV_X^D$--module. This constructions are used to prove the main theorem in this context. This theorem, proved by  Calderón-Moreno and Narváez-Macarro (\cite[Cor. 4.2]{calde_nar_fourier} and \cite[Corollary 1.7.2]{nar_comp_08}), states that a free divisor $D\subset X$ satisfies the Logarithmic Comparison Theorem if and only if the natural morphism
$$ \varrho:\DD_X \Lotimes_{\VV_X^D} \OO_X(D) \longrightarrow \OO_X(\star D)
$$ is an isomorphism in the derived category of $\cD_X$-modules (see Theorem \ref{th:D-crit-LCT}). Fixing a point $p\in D$, a reduced local equation $f$ of the germ $(D,p)$ and a basis $\{\delta_1,\ldots,\delta_n\}$ of $\fDer(-\log D)_p$, the morphism $\varrho_p$ is an isomorphism if and only if the complex
$\DD_{X,p} \Lotimes_{\VV_{X,p}^D} \OO_{X,p}(D)$ is concentrated in cohomological degree 0, $\cD_{X,p}f^{-1} = \cO_{X,p}(\star D)$  and the $\cD_{X,p}$--annihilator of $f^{-1}\in \cO_{X,p}(\star D)$ equals the left ideal $\cD_{X,p}(\delta_1+\alpha_1,\ldots,\delta_n+\alpha_n)$, where $\delta_i(f)=\alpha_i f$ for $i=1,\ldots, n$, see Sect. \ref{D-module-criterion-for-LCT}.  That gives a $\cD$--module criterion to test if LCT holds for a given free divisor $D\subset X$. Finally we treat some examples to show how this criterion can be applied in practice.

\section{Comparison Theorems}\label{compare}
\subsection{Grothendieck's Comparison Theorem\index{Grothendieck's Comparison Theorem}}\label{Grothendieck-comparison-theorem}
To calculate the cohomology of a space with constant coefficients $\CC$, one can use a resolution of the constant sheaf $\CC_X$ and then calculate the hypercohomology of the resolution. For example, on a complex manifold the holomorphic Poincar\'e lemma (see e.g. \cite[Ch. 2, Par. 2, Ex. 2.5.1]{Godement-Th-des-faisceaux} -- the proof for $C^\infty$ forms is easily adapted to the holomorphic case) shows that the complex of holomorphic differential forms
\beq\label{res1}\xymatrix{0\ar[r]&\CC_X\ar[r]&{\OO_X=\Omega}^0_X\ar[r]^d&{\Omega}^1_X\ar[r]^d&\cdots}\eeq
is exact, so that ${\Omega}^\bu_{X}$ is a resolution of $\CC_X$ by locally free $\OO_X$-modules. If $X$ is in addition a Stein space -- for example, the complement of a divisor $D$ in $\CC^n$, or a convex open set in $\CC^n$ -- then all of the sheaves in the complex are {\it acyclic}: $H^j(X,{\Omega}^k_X)=0$ for $j>0$ and $k\geq 0$ (in effect, this is the definition of `Stein space'). In this case the double complex with which one calculates the hypercohomology is reduced to the complex of global sections
$$0\tos\G(X,\OO_X)\tos \G(X,{\Omega}^1_X)\tos\cdots$$ and so one has the {\it analytic de Rham theorem}\index{Analytic de Rham theorem}:

\begin{teorema}If $X$ is a Stein manifold then \beq\label{iso1}\xymatrix{h^q(\G(X,{\Omega}^\bu_X))\ar[r]& H^q(X;\CC)}\eeq
is an isomorphism.
\end{teorema}

Here if we view $H^j(X;\CC)$ as the singular cohomology group, then the arrow is given by the integration of differential forms along
singular chains, which we will refer to as the {\it de Rham morphism}\index{de Rham morphism}. Note that this theorem already implies the non-trivial fact that the cohomology of a Stein manifold vanishes above middle dimension.
\vsn
If $X$ is the complement of a divisor (a hypersurface) $D\subset\CC^n$, then ${\Omega}^k_X$ (or more precisely $j_*{\Omega}^k_X$, where $j:X\tos\CC^n$ is inclusion) is equal to the sheaf of germs of holomorphic forms with arbitrary singularities along $D$. This is strictly larger than the sheaf ${\Omega}^k(\star D)$ of meromorphic forms with poles along $D$, since it places no restriction on the behaviour of the extension to $D$, and in particular allows also forms with essential singularities along $D$.  In \cite[Th. 2]{GRO}\footnote{It is noted in \cite{GRO} that the contents of
this paper formed part of a letter of the author to M. F. Atiyah, dated October 14, 1963, except for some remarks dated November 1963 and July 1965.}, Alexander Grothendieck showed that despite this difference, there are canonical isomorphism of sheaves of complex vector spaces
 \beq\label{iso3bis}
 h^q \left({\Omega}^\bu_{\CC^n}(\star D)\right)\simeq \RR^q j_*\CC_{\CC^n\ssm D},\quad q=0,\dots, n,\eeq where $\RR^q j_*$ is the $q$-th right derived functor of the left-exact functor $j_*$, see e.g. \cite[I.7]{iversen-1986}.
The maps in (\ref{iso3bis}) are given by the composition of the adjunction $${\Omega}^\bu_{\CC^n}(\star D) \tos  j_* j^{-1} {\Omega}^\bu_{\CC^n}(\star D) = j_*{\Omega}^\bu_{\CC^n\ssm D}$$ with the inverse of Poincaré quasi-isomorphism $\CC_{\CC^n\ssm D} \stackrel{\sim}{\tos} {\Omega}^\bu_{\CC^n\ssm D}$.

Taking into account that $\CC^n$ is Stein, by \cite[Corollary, page 97]{GRO}, one has  $H^j(\CC^n,{\Omega}^k_{\CC^n}(\star D))= 0 $ for $j>0$ and any $k$, and there are canonical isomorphisms of complex vector spaces
\beq\label{iso2}\xymatrix {h^q(\G(\CC^n,{\Omega}^\bu_{\CC^n}(\star D)))\ar[r]& H^q(\CC^n\ssm D;\CC)} \eeq for $q\geq 0$.

These maps are formal in the sense that they come from adjunction and Poincaré lemma. A well known (but not obvious) fact, see Remark \ref{comment-on-deRham-morphism}, is that these maps coincide with the de Rham map given by integration of meromorphic forms along singular chains composed with the isomorphism between singular cohomology with coefficients in $\CC$ and sheaf cohomology with coefficients in the constant sheaf $\CC_{\CC^n\ssm D}$ (see \cite[Ch. 4; Th. 4.14]{Ramanan2005}).

The isomorphism in (\ref{iso3bis}) can be stated as an isomorphism in the derived category, which is of course valid on any complex manifold $X$ and any divisor $D\subset X$, namely
\beq\label{iso3}{\Omega}^\bu_X(\star D)\simeq \RR j_*\CC_{X\ssm D}.\eeq
This is known as
{\it Grothendieck's Comparison Theorem} (see Theorem \ref{teo:Gro-comp-th-revisited}). The core assertion is that one can ignore essential singularities along $D$.

To end this section, let us emphasize that Grothendieck's proof of (\ref{iso3}) used in an essential way Hironaka's resolution of singularities (in the complex algebraic and the complex analytic settings), and became a model to follow for the cohomological study of algebraic varieties over an arbitrary field. This fact brought conferred upon the resolution of singularities in positive characteristic a privileged place in the construction of cohomological theories, to the extent that it appeared as a matter of substance in this context.
On the other hand, Grothendieck's comparison theorem became a completely new and unexpected way to understand regular singular points (at infinity) of integrable connections on algebraic fiber bundles (see footnote (13) in \cite{GRO}, dated July 1965, and the crucial work \cite{deligne-LNM}), and more generally, a conceptual way to incorporate {\em regularity} as a central notion in $\DD$-module theory (see \cite{meb_1980}).
Finally, $\DD$-module theory itself has provided new tools to understand Grothendieck's comparison theorem and to avoid resolution of singularities in its treatment \cite{meb_IHES, meb_1990} (see also \cite{meb_CIMPA_96}).

In this paper we would like: (1) to explain what the logarithmic comparison means, drawing some analogies with Grothendieck's comparison theorem; (2) to survey the original proof of the Logarithmic Comparison Theorem for locally quasihomogeneous free divisors (see Theorem \ref{lct}); and (3) to explain how $\DD$-module theory provides a characterization of those free divisors which satisfy the Logarithmic Comparison Theorem (see Sect. \ref{D-module-criterion-for-LCT}). We emphasise that in all of them resolution of singularities does not play any r\^ole.

\begin{nota}\label{comment-on-deRham-morphism}
For any topological space $X$, one denotes by $H_q(X;\CC)$ the singular homology group of $X$ with values in $\CC$, for $q \geq 0$. Recall that the singular homology group $H_q(X;\CC)$
is by definition the homology group $h_q(C_\bullet(X),\partial)$ of the complex $(C_\bullet(X),\partial)$ of singular chains of $X$ with coefficients in $\CC$, where $\partial$ denotes the boundary map.

Recall also that the complex of singular cochains $(C^\bullet(X),\delta)$ is defined as $$C^\bullet(X):=\Hom_\CC(C_\bullet(X), \CC),$$
where $\delta$ is the coboundary map associated with the boundary map $\partial$. The singular cohomology group $H^q(X;\CC)$, of $X$ with values in $\CC$ for $q\geq 0$, is by definition the cohomology group $h^q(C^\bullet(X),\delta)$.

There exists a natural morphism $$ H^q(X;\CC) \stackrel{\psi_q}{\longrightarrow} H_q(X;\CC)^*:= \Hom_\CC(H_q(X;\CC),\CC).$$ For any $[\eta]\in H^q(X;\CC)$ and any $[\gamma]\in H_q(X;\CC)$, one has $$\psi_q([\eta])([\gamma])= \eta(\gamma)$$ where $[\,]$ means equivalence class in the corresponding (co)homology group. The morphism $\psi_q$ is an isomorphism, for $q\geq 0$, since $\CC$ is a field.

If $X$ is a complex manifold with $\dim X =n\geq 1$,  the de Rham morphism\index{de Rham morphism} $$dR_q : h^q(\Gamma(X,\Omega_X^\bullet)) \longrightarrow H^q(X;\CC)$$ is defined (or induced) by integration of differential $q$-forms over singular $q$-chains: $$dR_q([\omega])([\gamma]) := \int_\gamma \omega \in \CC$$ where $\omega$ is a $q$-form, $\gamma$ is a singular $q$-chain and $[\,]$ means equivalence class in the corresponding (co)homology group. Here, to be precise, one needs to consider $C^\infty$ chains, and the point is that the inclusion of the complex of $C^\infty$ singular chains of $X$ with coefficients in $\CC$ in $C_\bullet(X)$ is a quasi-isomorphism (see \cite[Ap. A, Th. 2.1]{massey}).

The de Rham theorem\index{Analytic de Rham theorem} states that if $X$ is a Stein manifold then $dR_q$ is an isomorphism for $q=0,\ldots, n$; see e.g. \cite[App. A,  Th. 3.1]{massey}.
In the $C^\infty$ category, the de Rham theorem holds in much greater generality, for all paracompact differentiable manifolds, (see {\it loc. cit.}).

Recall that by the Poincaré Lemma the complex (\ref{res1}), i.e. the complex $$ 0 \rightarrow \CC_X \rightarrow \Omega_X^\bullet,$$  is  exact. Since the category $\Mod(\CC_X)$, of sheaves of $\CC$-vector spaces, has enough injectives (see e.g. \cite[Ex. 2.3.12]{weibel-1994}), there exists an injective resolution $$0 \rightarrow  \CC_X \rightarrow I_X^\bullet$$ and a natural morphism $$ \Omega_X^\bullet \rightarrow I_X^\bullet,$$ unique up to homotopy equivalence, making the obvious diagram commute (see \cite[Rq. 3, Th. 2.4.1]{GRO-Tohoku} and e.g. \cite[2.5]{weibel-1994}). So for any $q$,  there is a natural morphism
$$\alpha_q: \,  h^q(\Gamma(X,\Omega_X^\bullet)) \rightarrow h^q(\Gamma(X,I^\bullet_X)) = H^q(X;\CC_X)$$ where the last equality is just the definition of the sheaf cohomology of the constant sheaf $\CC_X$.

Moreover, as mentioned in the introduction, if $X$ Stein, then each $\Omega^q_X$ is acyclic for the ``global sections functor'' $\Gamma(X;\,)$, and hence the natural morphism $\alpha_q$ is an isomorphism (see \cite[Rq. 3, Th. 2.4.1]{GRO-Tohoku} and e.g. \cite[2.5]{weibel-1994}).

Finally, for any complex manifold $X$ (in fact, for any locally contractible topological space) let us consider, see e.g. \cite[Ch. 4, Th. 4.14]{Ramanan2005}, the complex
\begin{equation}\label{singcomp}0 \rightarrow \CC_X \rightarrow \widetilde{S}_X^\bullet\end{equation}
where, for $q\geq 0$, one denotes by $\widetilde{S}_X^q$ the sheaf associated to the presheaf ${S}_X^q$ of singular $q$--cochains on $X$. Local contractibility means that each point has a fundamental system of neighbourhoods which are contractible. In each of these, the complex of singular cochains is exact. It follows that the complex of sheaves \eqref{singcomp} is exact. As before for the complex $\Omega^\bullet_X$ (but now with the complex $\widetilde{S}_X^\bullet$ in its place), and again using the injective resolution $I^\bullet_X$ of $\CC_X$,   for each $q$ there is a natural morphism
$$\beta_q: \, H^q(X;\CC) = h^q(\Gamma(X,\widetilde{S}_X^\bullet)) \rightarrow h^q(\Gamma(X,I^\bullet_X)) = H^q(X;\CC_X)$$ where the first equality is just the definition of the singular cohomology $H^q(X;\CC)$ (see \cite[Ch. 4, Prop. 4.12]{Ramanan2005}). Since the sheaves $\widetilde{S}_X^q$ are all flabby they are also acyclic for the functor $\Gamma(X;\,)$ (see \cite[Ch. 4, ex. 1.10, 2 and Prop. 3.3]{Ramanan2005}). Hence the morphism $\beta_q$ is an isomorphism for  $q\geq 0$.

So, for any $q$ we have the following diagram of isomorphisms
\beq \label{eq:triangle-de-Rham}
\begin{tikzcd}
h^q(\Omega^\bullet (X))  \ar[]{r}{\text{$dR_q$}} \ar[]{d}{\text{$\alpha_q$}} &  H^q(X;\CC) \ar[]{dl}{\text{$\beta_q$}} \\ 
H^q(X;\CC_X).  &  
\end{tikzcd}
\eeq
It is ``well-known'' that the diagram commutes, up to a sign, although a concrete reference for this result seems to be elusive.
\end{nota}

\subsection{The Brieskorn Complex\index{Brieskorn Complex}}\label{sub:Brieskorn-complex}

For a special class of divisors, namely hyperplane arrangements in affine space $X=\CC^n$, Brieskorn, in \cite{brie} generalising Arnol'd in \cite{ar69}, had already given a way of calculating the cohomology of the complement. If $D$ is the union of hyperplanes $H_j$ with equations $h_j=0$, $j=1,\ld, N$, then the collection of meromorphic 1-forms $\frac{dh_j}{h_j}$ generates a $\CC$-subalgebra $B^\bu$ of the exterior algebra $\Gamma(X,\Omega_X^\bullet(\star D))$. Since $d\left(\frac{dh_j}{h_j}\right)=0$, all exterior derivatives on $B^\bu$ are zero, and $B^\bu$ is a subcomplex of the complex $\Gamma(X,\Omega_X^\bullet(\star D))$ with derivative zero, known as the {\it Brieskorn complex}. Brieskorn showed

\begin{teorema} The de Rham morphism $B^p\tos H^p(U;\CC)$ is an isomorphism.\hspace*{\fill}$\Box$ \vskip \baselineskip
\end{teorema}

Note that each of Brieskorn's forms $$\omega_{j_1,\ld,j_k}:=\frac{dh_{j_1}}{h_{j_1}}\wedge\cdots\wedge \frac{dh_{j_k}}{h_{j_k}}$$ has at most a first order pole along $D$, so in this special case his result is stronger than Grothendieck's.

An earlier version of Brieskorn's result, that it holds for a normal crossing divisor (NCD)(\cite{AH55})  played an important role in Deligne's mixed Hodge theory: every quasiprojective variety $U$ has a {projective compactification}  $U\hookrightarrow X$ in which $D:=X\ssm U$ is a NCD. Similarly, every singular variety has a resolution whose exceptional divisor is once again an NCD. Deligne used  the poles to define the weight filtration on the cohomology of $U$.

If $D=\bigcup_{j=1}^N\{h_j=0\}$ is a union of more general irreducible divisors $\{h_j=0\}$, the analogous Brieskorn complex $B^\bu$ generated by the closed forms $\frac{dh_j}{h_j}$ does {\it not} in general calculate the cohomology of the complement: for example, if $D$ is irreducible, this complex reduces to $0\tos B^0\tos B^1\tos 0$, whereas $X\ssm D$ may have cohomology in higher dimensions. Deligne and Dimca
explored the relation between pole order and the Hodge filtration in \cite{deldim}.

\subsection{The Logarithmic Comparison Theorem\index{Logarithmic Comparison Theorem}}\label{LCTfirst}
Kyoji Saito, in \cite{Saito-80},  introduced a fruitful generalisation of the Brieskorn complex. A meromorphic form $\omega$ has a {\it logarithmic pole} along the divisor $D$ if both $h{\omega}$ and $hd{\omega}$ are regular, where $h$ is a local equation for $D$. Note in particular that all of Brieskorn's forms have logarithmic poles along $D$.

Denote the sheaf of germs of meromorphic $k$-forms with logarithmic poles\index{Logarithmic differential form} by $\Omega^k_X(\log D)$. Then $\Omega_X^\bu(\log D)$ is again a subcomplex of  $\Omega_X^\bu(\star D)$. In the same article Saito also described the dual of $\Omega_X^1(\log D)$, namely the sheaf of \,``logarithmic derivations''\index{Logarithmic derivation},\,  $\fDer(-\log D)$, whose stalk at $x\in D$ consists of germs of vector fields on $X$ which are tangent to $D$ at its smooth points. The duality arises from the contraction pairing
$$\Omega_X^1(\log D)\times \fDer(-\log D)\tos \OO_X, \quad({\omega},\xi)\mapsto \iota_\xi(\omega)={\omega}(\xi).$$
Saito showed that this is a perfect pairing, so that $\Omega_X^1(\log D)$ is also the dual of $\fDer(-\log D)$.  Saito's interest was especially in the case where $D$ is a {\it free divisor}\index{Free divisor}, that is, where $\fDer(-\log D)$ (or, equivalently, $\Omega_X^1(\log D)$) is a locally free sheaf of $\OO_X$-modules. Because the dual of any $\OO_X$-module has depth at least 2 when $\depth\, \cO_X \geq 2$ (see e.g. \cite{schless}, or
\cite[Lemma 9.2]{looij}), it follows from this that every plane curve is a free divisor. More interesting is the fact, proved by K. Saito, that the discriminant in the base of a versal deformation of an isolated hypersurface singularity is a free divisor. Looijenga generalised this by showing in \cite[Corollary 6.13]{looij} that it holds also for the discriminant in the base of a versal deformation of an isolated complete intersection singularity (ICIS), and later authors (\cite{vans}, \cite{bebo}) have extended this to the discriminants of a range of non-ICIS singularities, to discriminants in quiver representation spaces \cite{bm} and more generally in prehomogeneous vector spaces \cite{gms}.

Since for any divisor $D\subset X$ one has a natural morphism $$ \Omega_X^\bullet(\log D) \hookrightarrow \Omega_X^\bullet(\star D)$$ one can ask for the class of divisors $D\subset X$ such that the previous natural morphism is a quasi-isomorphism (or an isomorphism in the derived category of sheaves of $\CC_X$--vector spaces).

\begin{definicion} \label{def:LCT-holds-for-D}
If the morphism $ \Omega_X^\bullet(\log D) \hookrightarrow \Omega_X^\bullet(\star D)$  is a quasi-isomorphism we say that the divisor $D$ satisfies the {\em logarithmic comparison theorem} (LCT), or even that {\em the LCT holds for $D$}.
\end{definicion}

By Grothendieck's comparison theorem, the LCT holds for a divisor $D$ if and only if the morphism $h^q(\Omega_X^\bullet(\log  D)) \rightarrow \RR^q j_* \CC_{U}$ is an isomorphism for $0\leq q \leq \dim X$.
Any NCD satisfies the LCT, after \cite[Cap. II, Lemme 6.9]{deligne-LNM}. The result of Atiyah and Hodge cited above shows that every cohomology class on the complement of the divisor can be represented by a logarithmic form, but not that the logarithmic complex calculates this cohomology.

\begin{teorema}\label{LCT}{\rm [Logarithmic Comparison Theorem\index{Logarithmic Comparison Theorem}, \cite{CMN}]}\label{lct}  If $D\subset X$ is a locally quasihomogeneous free divisor with complement $U$, then the de Rham morphism
\beq\label{iso4}{\Omega}^\bu_X(\log D)\tos \RR j_*\CC_U\eeq
is a quasi-isomorphism.
\end{teorema}

Here {\it locally quasihomogeneous}\index{Locally quasihomogeneous} means that at each point $p\in D$ there are local coordinates on $X$ in which $D$ has a weighted homogeneous local equation with all weights strictly positive -- in other words, locally there is a good $\CC^*$-action centred at $p$.   Divisors with this property are also known as {\it strongly quasihomogeneous}, and as {\it positive}.  Evidently hyperplane arrangements have this property, and it also holds for the discriminants of stable maps in Mather's ``nice dimensions'', because of the remarkable (and unexplained) fact that in the nice dimensions, all stable germs are quasihomogeneous in suitable coordinates (see \cite[Section 7.4]{mn} for a pedestrian proof of this).
\vsn
Local quasihomogeneity is used twice in the proof of Theorem \ref{lct}. First, it allows an inductive argument on the codimension of strata.

\begin{lema}\label{sqh}
Let $X$ be a
complex manifold of dimension $n$, let $D$ be a
strongly quasihomogeneous divisor in $X$, and let $p\in D$. Then there is an
open neighbourhood $U$ of $p$ such that for each $q\in U \cap D$, with $q\neq  p$, the germ at $q$ of the
pair $(X,D)$ is isomorphic to the germ at $0$ of a product
$(\CC^{n-1}\times\CC, D_0\times\CC)$ where $D_0$ is a strongly quasihomogeneous divisor.
\end{lema}
\vsn
To prove that \eqref{iso4} is an isomorphism at $p\in D$, one uses induction on the dimension of the divisor. It is easy to see that if \eqref{iso4} holds for a a divisor $D_0\subset \CC^{n-1}$ then it holds for $D_0\times\CC\subset \CC^{n-1}\times\CC$. Lemma \ref{sqh} therefore says that we may assume by induction that \eqref{iso4} holds at all points of $D\cap U\ssm\{p\}$. Note that the induction begins with the divisor
$0\subset\CC$; it is well known that $H^1(\CC\ssm\{0\};\CC)\simeq\CC$, generated by the logarithmic form $\frac{dz}{z}$. Thus in what follows we assume $n\geq 2$.
\vsn
The main body of the proof of \ref{lct} then consists in showing that \eqref{iso4} also holds at $p$ (which for convenience we suppose is the point $(0,\ld,0)$), by showing that for any sufficiently small polycylinder $V$ centred at $p$, the de Rham morphism $$h^p(\Gamma(V,\Omega_X^\bu(\log D))\tos H^p(V\ssm D;\CC)$$ is an isomorphism. The argument involves two \v Cech-de Rham double complexes associated with Stein covers
$\{V_i\}$ of $V\ssm\{0\}$, with $V_i=V\ssm\{x_i=0\}$, and $\{V_i'=V_i\ssm D\}$ of $V\ssm D$, and the two standard spectral sequences associated to each one. The complexes are
\begin{align}K^{p,q}&=\bigoplus_{i_0<i_1<\cdots<i_q}\Gamma\left(\bigcap_{j=0}^qV_{i_j}, {\Omega}_X^p(\log D)\right)\label{k}\\
\intertext{and}
\tilde K^{p,q}&=\bigoplus_{i_0<i_1<\cdots<i_q}\Gamma\left(\bigcap_{j=0}^qV'_{i_j}, {\Omega}_X^p\right)\label{tildek}
\end{align}
with differentials $d$, the exterior derivative, and $\check d$, the Cech differential. The inclusion $V_i'\hookrightarrow V_i$ determines a morphism $\rho_0:K^{p,q}\tos \tilde K^{p,q}$ which commutes with $d$ and $\check d$ and thus gives rise to morphisms of the associated spectral sequences, which we denote by $'\rho_\ell$ and $''\rho_\ell$, where the subindex $\ell$ refers to the page of the spectral sequence.
\vsn
{\bf First spectral sequences}
\vsn
Applying $d$ to \eqref{k} and \eqref{tildek}, we get the first page of the first spectral sequence associated to each double complex:
\begin{align*}'E_1^{p,q}&=h^p(K^{\bu,q})=\bigoplus_{i_0<i_1<\cdots<i_q}h^p\left(\Gamma\left(\bigcap_{j=0}^qV_{i_j}, {\Omega}_X^\bu(\log D)\right)\right)
\\
\intertext{and}
'\tilde{E}_1^{p,q}&=h^p(\tilde K^{\bu,q})=\bigoplus_{i_0<i_1<\cdots<i_q}h^p\left(\Gamma\left(\bigcap_{j=0}^qV'_{i_j}, {\Omega}_X^\bu\right)\right)=
\bigoplus_{i_0<i_1<\cdots<i_q}H^p(\bigcap_{j=0}^qV'_{i_j};\CC).
\end{align*}
Since $0$ is excluded from all of the open sets in the covers, the induction hypothesis implies that $'\rho_1^{p,q}:{'E}_1^{p,q}\tos {'\tilde E}_1^{p,q}$ is an isomorphism for all $p,q$, and it follows that
$'\rho_\infty^{p,q}$ is also an isomorphism. Thus, $\rho_0$ induces an isomorphism of the cohomology of the total complexes of
$K^{\bu,\bu}$ and $\tilde K^{\bu,\bu}$.
\vsn
{\bf Second spectral sequences}
\vsn
Because $\{V_i\}$ and $\{V_i'\}$ are Stein covers, applying $\check d$ gives,  as first pages of the second spectral sequences,
$$
''E^{p,q}_1=H^q(V\ssm\{0\}, {\Omega}_X^p(\log D))\quad \text{and}\quad
''\tilde E^{p,q}_1=H^q(V\ssm D, {\Omega}_X^p).
$$
Because $D$ is a free divisor, all of the ${\Omega}_X^p(\log D)$, as well as the ${\Omega}_X^p$, are free $\OO_X$-modules, and thus since
$H^q(V\ssm\{0\},\OO_X)$ is zero except for $q=0$ and $q=n-1$ (see e.g. \cite[(8.14)]{looij}), $H^q(V\ssm\{0\},{\Omega}_X^p(\log D))$
also vanishes except in these dimensions. In other words
$$''\tilde E_1^{p,q}=0 \text{\ \ except for } q=0\text{ and }q=n-1.$$
Now because $V\ssm D$ is a Stein space and ${\Omega}_X^p$ is coherent, $''\tilde E^{p,q}_1=H^q(V\ssm D,{\Omega}_X^p)$ is equal to $0$ for all $q>0$, and the spectral sequence $''\tilde E$ collapses at $E_2$, with
\beq\label{niso''}''\tilde E_\infty^{p,0}=H^p(V\ssm D;\CC),\ \ ''\tilde E^{p,q}_\infty=0\ \text{ if }q>0.\eeq
We claim that $''E$ also collapses at $E_2$. In view of the vanishing already remarked upon, it is enough to show that the complex
$(''E_1^{\bu,q},d)$, i.e.
\beq\label{les}0\tos H^{n-1}(V\ssm\{0\},{\Omega}_X^0(\log D))\tos H^{n-1}(V\ssm\{0\},{\Omega}_X^1(\log D))\tos\cdots\tos H^{n-1}(V\ssm\{0\},{\Omega}_X^t(\log D))\tos 0\eeq is exact.
This is the only point at which the argument departs from ``general nonsense''.

Each of the groups in \eqref{les} is generated by
classes of the form $c_{\alpha}{\omega}_{\alpha}$ where $c_{\alpha}\in H^{n-1}(V\ssm\{0\},\OO_X)$ and ${\omega}_{\alpha}$ is a basis element of ${\Omega}_X^p(\log D)$ in $V$.

Recall that $(D,0)$ is assumed weighted homogeneous in $V$, with respect to positive weights $w_1,\ld, w_t$ for the coordinate functions.
We have $$H^{n-1}(V\ssm \{0\},\OO_X)=\frac{\Gamma(V\ssm \bigcup_i\{x_i=0\},\OO_X)}{
\sum_j\Gamma(V\ssm \bigcup_{i\neq j}\{x_i=0\},\OO_X)}.$$
Each term in the numerator can be represented by a Laurent series in which all exponents are negative; series in the denominator have all exponents negative, bar one.  An easy lemma shows that ${\Omega}_X^p(\log D)_0$ has a basis consisting of forms of weighted degree strictly less than
$\sum_jw_j$. It follows that {\it the complex \eqref{les} has zero weight-zero part}. But since the exterior derivative $d$ preserves weighted degrees, it follows that using the Lie derivative with respect to the Euler form, one can construct a contracting homotopy from the complex to its weight zero part. Thus, it is exact.

We are left with
\beq\label{iso''}''E_\infty^{p,0}=h^p(\Gamma(V\ssm\{0\}, \Omega_X^\bu(\log D))=h^p(\Gamma(V, \Omega_X^\bu(\log D)), \ \ ''E_\infty^{p,q}=0\ \text{ if }q>0\eeq
(since $tn>1$ and $\Omega^p(\log D)$ is locally free, $0$ is a removable singularity).
The theorem now follows  from \eqref{niso''}, \eqref{iso''} and the isomorphism of the cohomology of the total complexes of $K^{\bu,\bu}$ and $\tilde K^{\bu,\bu}$.

\vsn
The importance of local quasihomogeneity for logarithmic comparison theorem remains unclear. It was shown to be necessary for plane curves in \cite{cmnc}, and for free surfaces in 3-space in \cite{GrangerSchulze}, but in higher dimensions, local weak quasihomogeneity is sufficient, see \cite{wqh}, \cite[Remark 1.7.4]{nar_comp_08} {and  Theorem \ref{th:LCT-LWQH}}.

\begin{definicion} \label{def:LWQH}
A divisor $D$ is locally weakly quasihomogeneous\index{Locally weakly quasihomogeneous divisor} \cite[Def. 2.1]{wqh} if at each point $p\in D$ there are local coordinates on $X$ in which $D$ has a weak weighted homogeneous local reduced equation $f=0$, that is: all the weights are non negative and not all of them are 0, with $f$ of strictly positive weight.
\end{definicion}

\subsection{Conjecture of Terao\index{Conjecture of Terao}}\label{conjecture-of-terao} H. Terao conjectured in \cite[Conjecture 3.1]{terao} that LCT always holds for hyperplane arrangements.  After partial results of Wiens and Yuzvinsky in \cite{wy}, the conjecture was finally proved by Daniel Bath in \cite{Bath2022}. Bath's proof uses induction on the codimension of the flats of the arrangement, and a spectral sequence argument similar to the argument of \cite{CMN} given above.
\subsection{Linear free divisors and logarithmic comparison revisited\index{Linear free divisors}}\label{linear-free-divisors}

\noindent A free divisor $D\subset \CC^n$ is {\it linear} if $\Der(-\log D) = \Gamma(\CC^n,\fDer(-\log D))$ has a basis of vector fields whose coefficients are linear forms. Normal crossing divisors are of course the simplest examples, and are in fact the only linear free divisors among hyperplane arrangements. In \cite{bm} R.-O. Buchweitz and the second author showed that the discriminant in the representation space of a Dynkin quiver, with a root of the underlying diagram as dimension vector, is a linear free divisor. All irreducible linear free divisors are classified, though not under that name, in \cite{sk}, where the operation of {\it castling} derives them from four basic examples.

If $D\subset \CC^n$ is a linear free divisor then
the set of weight zero vector fields in $\Der(-\log D)$ is a Lie algebra of dimension $n$; it is isomorphic to the Lie algebra of the Lie subgroup of $G_D$ consisting of
linear automorphisms of the pair $(\CC^n,D)$. For linear free divisors, \cite{gmns} gives an argument for the logarithmic comparison theorem which is quite different from the proof of \cite{CMN}.  When D is a linear free divisor, then it turns out that $\CC^n\ssm D$ is a single orbit of the identity component $G^0_D$ of the group $G_D$, so $H^*(\CC^n\ssm D;\CC)$ is the cohomology of $G^0_D/S_p$, where $S_p$ is the isotropy subgroup of a point $p$. Since $G^0_D$ is path connected, the action of $S_p$ on $H^*(G^0_D;\CC)$ is trivial, so $H^*(G^0_D/S_p;D)=H^*(G^0_D;\CC)$. Second, in this case the weight zero subcomplex of $\G(\CC^n,{\Omega}_{\CC^n}^\bu(\log D))$ coincides with the complex of Lie algebra cohomology, with complex coefficients, of the Lie algebra ${\mathfrak{g}}_D$ of $G^0_D$. By an argument using weighted homogeneity, $\G(\CC^n,{\Omega}_{\CC^n}^\bu(\log D))$ is chain homotopic to its weight zero subcomplex. Thus, $\G(\CC^n,{\Omega}_{\CC^n}^\bu(\log D))$ coincides with the Lie algebra cohomology of ${\mathfrak{g}}_D$ with complex coefficients. For compact connected Lie groups $G$, a well-known argument shows that the Lie algebra cohomology coincides with the topological cohomology of the group. For linear free divisors the group $G^0_D$ is never compact, but the isomorphism also holds good for the larger class of reductive groups, and for a significant class of linear free divisors, including all those mentioned above, $G^0_D$ is indeed reductive.
\subsection{A pairing $\mathbf{H^p(X\ssm D)\times H^q(D)\tos H^{p+q}(D)}$} \label{sec:omega-cech}

The paper \cite{dfafd} introduces a variant $\check{\Omega}^k_D$ of the module of K\"ahler forms ${\Omega}^k_D$, defined by
$$\check{\Omega}^k_D:=\frac{\Omega^k_X}{h\Omega_X^k(\log D)},$$
where $h$ is a reduced equation for $D$.
Note that $\check{\Omega}^0_D={\Omega}^0_D=\OO_D$.
Since $\frac{dh}{h}$ has a logarithmic pole and $ {\Omega}_X^k(\log D)\supset {\Omega}^k_X$, it follows that $h{\Omega}_X^k(\log D)\supseteq dh\wedge{\Omega}^{k-1}_X+h{\Omega}^k_X$, so $\check{\Omega}^k_D$ is a quotient of ${\Omega}^k_D$. If $D$ is a free divisor then each $\check{\Omega}^k_D$ is a maximal Cohen-Macaulay $\OO_D$-module, and coincides with ${\Omega}^k_D$ at smooth points of $D$. The $\check{\Omega}^k_D$ form a complex with respect to the usual exterior derivative:
if $\omega\in{\Omega}_X^k(\log D)_x$ then $\frac{dh}{h}\wedge\omega\in{\Omega}_X^{k+1}(\log D)$, from which it follows that $d(h\omega)\in h\Omega_X^{k+1}(\log D)$.  We denote the exterior derivative on this complex by $\check d$.
If $D$ is quasihomogeneous then $(\check{\Omega}^\bu_D, \check d)$ is a resolution of $\CC_D$ (\cite[Lemma 3.3]{dfafd}).
\vsn
Straightforward calculations show
\begin{enumerate}
\item
there is a well-defined
pairing \beq\label{pairing}{\Omega}_X^k(\log D)\times \check{\Omega}^\ell_D\tos \check {\Omega}^{k+\ell}_D,\ \ \text{defined by }
({\omega}_k, {\omega}_\ell)\mapsto {\omega}_k\wedge{\omega}_\ell.\eeq
Note that there is no comparable wedge pairing
$${\Omega}_X^k(\log D)\times {\Omega}^\ell_D\tos {\Omega}^{k+\ell}_D.$$
For  such a pairing to be well defined, the wedge of ${\omega}_k\in {\Omega}_X^k(\log D)$ with $\omega_\ell\in(dh\wedge{\Omega}^{\ell-1}_X+h{\Omega}^\ell_X)$ would have to lie in $dh\wedge{\Omega}^{k+\ell-1}_X+h{\Omega}^{k+\ell}_X$, and in general this does not hold.
For example, if $h=h_1h_2$, and we take $\omega_1=\frac{dh_1}{h_1}\in\Omega_X^1(\log D)$ and ${\omega}_1'=dh\in dh\wedge{\Omega}^0_X$, then
$${\omega}_1\wedge{\omega}_1'=\frac{dh_1}{h_1}\wedge \bigl(h_1dh_2+h_2dh_1\bigr)=dh_1\wedge dh_2\notin \bigl(dh\wedge{\Omega}^1_X+h{\Omega}^2_X\bigr).$$
\item
The pairing \eqref{pairing} descends to a pairing on the homology of the two complexes,
$$h^k(\Gamma(X,{\Omega}_X^\bu(\log D))\times h^\ell(\Gamma(X,\check{\Omega}^\bu_D))\tos h^{k+\ell}(\Gamma(X,\check{\Omega}^\bu_D)).$$
\item
When $D$ is locally quasihomogeneous then in view of the LCT and the exactness of $\check {\Omega}^\bu_D$, this gives a pairing
$$H^k(X\ssm D;\CC)\times H^\ell(D;\CC)\tos H^{k+\ell}(D;\CC).$$
The properties of this pairing remain to be explored.
\end{enumerate}

\section{$\DD_X$--modules\index{$\DD$--modules}} \label{D_X-modules} In this section we recall some of the basics of $\DD$--module theory.
We mainly follow \cite{sem-gre-1975}, \cite{meb_formalisme} and \cite{Granger-Maisonobe-cimpa-niza}. The section contains the necessary terminology and results that are used in the subsequent sections.

\subsection{Basic objects} \label{basic-objects} Let $X$ be a complex manifold of dimension $n\geq 1$ and $\cO_X$ (or simply $\cO$) be the sheaf of germs of holomorphic functions on $X$. The sheaf of rings (or more precisely of $\CC$--algebras) $\cO_X$ is coherent: this is Oka's theorem \cite{Oka-Kiyoshi}, \cite{sem-cartan-51-52}.

We denote by $\DD_X$ (or simply by $\DD$) the sheaf of linear differential operators\index{Linear differential operator} on $X$ with holomorphic coefficients \cite[I, \S 1]{sem-gre-1975}, \cite[Def. 3]{Granger-Maisonobe-cimpa-niza}. Sometimes we say {\em operators (or differential operators) on $X$} instead of linear differential operators on $X$ with holomorphic coefficients. For any open set $U\subseteq X$, $\DD(U)$ is a  $\cO(U)$--algebra and so $\DD$ is a sheaf of $\cO$--algebras.

If $(U; x_1,\ldots,x_n)$ is a chart in $X$ and $U\subset X$ is connected, any operator $P\in \DD(U)$ can be written in a unique way as a finite sum $$P = \sum_{\alpha \in \NN^n} c_\alpha \partial^\alpha =  \sum_{\alpha \in \NN^n} c_\alpha \partial_1^{\alpha_1} \cdots  \partial_n^{\alpha^n}$$ where $\alpha = (\alpha_1,\ldots,\alpha_n)$, $c_\alpha \in \cO(U)$ and $\partial_i$ is the partial derivative $\frac{\partial}{\partial x_i}$. For any $x\in X$, the stalk $\DD_{X,x}$ is a non-commutative $\cO_{X,x}$--algebra, since $\partial_i x_i - x_i \partial_i = 1$ for $1\leq i \leq n$.

If $P\in \DD(U)$ is non zero, the order of $P$ is the non negative integer $$\ord(P)=\max \{|\alpha|\, ;\, c_\alpha \not= 0\}$$ where $|\alpha| = \alpha_1 +\cdots+\alpha_n$. We also write $\ord(0)=-\infty$.

For $k\in \NN$, we denote by $F^k(\DD)$ the subsheaf of $\DD$ whose sections can be written locally as operators of order less than or equal to $k$. Each $F^k(\DD)$ is a sheaf of coherent $\cO$--modules and the family $F^\bullet := (F^k(\DD))_k$ is a discrete increasing exhaustive filtration of $\DD$. One has $F^0(\DD)=\cO$ and $F^k(\DD)F^{\ell}(\DD)=F^{k+\ell}(\DD)$ for $k,\ell\geq 0$. The associated sheaf of graded rings (and more precisely, of graded $\CC$--algebras) is $$\gr_{F^\bullet}(\DD) : = \bigoplus_{k\in \NN} \frac{F^k(\DD)}{F^{k-1}(\DD)}$$ where we write $F^{-1}(\DD)=\{0\}$.

If $P,Q$ are local sections in $\DD$, the difference $[P,Q]:=PQ-QP$ is called the commutator of $(P,Q)$.
For $x\in X$, the stalk $\gr_{F^\bullet}(\DD_x) \simeq \gr_{F^\bullet}(\DD)_x$
is a commutative ring since $\ord([P,Q]) \leq \ord(P)+\ord(Q)-1$. Thus, $\gr_{F^\bullet}(\DD)$ is a sheaf of commutative $\CC$--algebras.

The sheaves of rings $\DD_X$ and $\gr_{F^\bullet}(\DD_X)$ are coherent (see \cite[I, Th. 3.2]{sem-gre-1975}, \cite[Prop. 9]{Granger-Maisonobe-cimpa-niza}).

We denote  $$\sigma_k : F^{k}(\DD) \rightarrow  \frac{F^k(\DD)}{F^{k-1}(\DD)}$$ the natural projection map (which is a morphism of $\OO$--modules) and we call it the $k^{th}$ {\em symbol map}.

If $P$ is an operator of order $k$ we simply write $$\sigma(P)=\sigma_k(P)$$ and we call this element {\em the principal symbol}\index{Principal symbol} of $P$.

If $(x_1,\ldots,x_n)$ is a system of local coordinates around a point  $x\in X$, then $\partial_i=\frac{\partial}{\partial x_i} \in F^1(\DD_x)$ and we write $\xi_i = \sigma(\partial_i)$ for $1\leq i \leq n$.

There is an isomorphism of $\cO_x$--algebras
\beq\label{grD-O-xi}
\gr_{F^\bullet}(\DD_x) \longrightarrow \cO_x[\xi_1,\ldots,\xi_n]
\eeq
see e.g. \cite[I, \S 1]{sem-gre-1975}, \cite[Prop. 4]{Granger-Maisonobe-cimpa-niza};  the image of an element $P+F^{k-1}(\DD_x)$, with $P\in F^k(\DD_x)$, being simply $\sigma_k(P)$. Compare this morphism with (\ref{symmetric-algebra}). Since $\gr_{F^\bullet}(\cD_x)$ is noetherian, a standard argument using induction on the order of the operators, proves that $\cD_x$ is a left and right noetherian ring.

\subsection{Filtrations\index{Filtration}}

Let $\cM$ be a left $\DD$--module (i.e. a sheaf of left $\DD$--modules). A filtration of $\cM$ is a collection, $(\cM_k)_{k\in \NN}$, of $\cO$-submodules of $\cM$ such that:
\begin{enumerate}
\item[(i)] For any $k,\ell \in \NN$ : $\cM_k \subset \cM_{k+1}$ and $F^{\ell}(\DD)\cM_k \subset \cM_{k+\ell}$ .
\item[(ii)] $\cM = \cup_{k\in\NN} \cM_k$.
\end{enumerate}

To any filtration $\Gamma:=(\cM_k)_{k\in \NN}$ of $\cM$ we associate its graded $\cO$--module $$\gr_{\Gamma}(\cM):= \bigoplus_{k\in \NN} \frac{\cM_k}{\cM_{k-1}}$$ with $\cM_{-1}:=\{0\}$. More generally, $\gr_{\Gamma}(\cM)$ is a $\gr_{F^\bullet}(\DD)$--module.
A filtration $(\cM_k)_k$ of $\cM$ is said to be {\em good} is the following two conditions hold
\begin{enumerate}
  \item For any $k\in\NN$, $\cM_k$ is a coherent $\cO_X$--module.
  \item There exist $k_0\in\NN$ such that for any $\ell\in \NN$, $F^\ell(\DD) \cM_{k_0} = \cM_{k_0+\ell}$.
\end{enumerate}

Any coherent $\DD$--module $\cM$ admits locally (i.e. after restriction to sufficiently small open subsets in $X$)  good filtrations (see e.g. \cite[I, p. 7]{sem-gre-1975}, \cite[Prop. 10]{Granger-Maisonobe-cimpa-niza}).

Let $\Gamma:=(\cM_k)_k)$ be a filtration of $\cM$. Then $\Gamma$ is locally good (i.e. after restriction to sufficiently small open subsets in $X$, $\Gamma$ is a good filtration) if and only if $\gr_{\Gamma}(\cM) $ is a coherent $\gr_{F^\bullet}(\DD)$--module (see \cite[I, Prop. 4.1]{sem-gre-1975}, \cite[Th. 1]{Granger-Maisonobe-cimpa-niza}).

Let $\cM$ be a coherent left $\cD_X$--module. Then there exists a coherent sheaf of ideals $\cJ(\cM)$ of $\gr_{F^\bullet}(\cD_X)$ with the following property: for any open subset $U\subset X$ such that the restriction $\cM_U$ admits a good filtration, one has $\cJ(\cM)_U \simeq \sqrt{\ann_{\gr(\cD_U)}(\gr(\cM_U))}$, where $\gr(\cD_U) = \gr_{F^\bullet}(\cD_U)$ and $\ann_{\gr(\cD_U)}(\gr(\cM_U))$ is the sheaf of annihilating ideals of the coherent $\gr(\cD_U)$--module $\gr(\cM_U)$ (see \cite[Prop.17; Rmk. 6]{Granger-Maisonobe-cimpa-niza}).

\begin{definicion} \index{Characteristic variety} \label{variedad-caracteristica} Let $\cM$ be a coherent $\cD_X$--module. The closed analytic subset of the cotangent bundle $T^*X$ defined by the coherent sheaf of ideals $\cJ(\cM)$ is called {\em the characteristic variety} of $\cM$ and is denoted by $\Char(\cM)$.
\end{definicion}

\begin{ejemplo} 1) The characteristic variety of the $\cD_X$--module $\{0\}$ is $\emptyset$, and $\Char(\cD_X)=T^*X$ and then $\dim \Char(\cD_X)=2n=2\dim X$. \\
2) If $\cM=\frac{\cD_X}{\cI}$ where $\cI \subset \cD_X$ is a coherent sheaf of left ideals in $\cD_X$, the characteristic variety $\Char(\cM)$ is the closed analytic subset of $T^*X$ defined by the coherent sheaf of ideals $\gr_{F^{\bullet}}(\cI) \subset \gr_{F^\bullet}(\cD_X)$, since the filtration on the quotient $\cM$, induced by the $F^\bullet$ filtration on $\cD_X$, is a locally good filtration.
\\ 3) The $\cD_{X}$--module $\cO_{X}$ can be presented locally as a quotient $\cD_{X,x}/\cD_{X,x}(\partial_1,\cdots,\partial_n)$ where $\{\partial_1,\ldots,\partial_n\}$ is a basis of $\fDer(\cO_{X,x})$. The graded ideal $\gr(\cD_{X,x}(\partial_1,\cdots,\partial_n))$ equals $\gr(\cD)(\xi_1,\ldots,\xi_n)$ (recall that $\xi_i$ is the principal symbol of $\partial_i$ (\ref{basic-objects})).  Then $\Char(\cO_X)=T^*_X X$, that is, the zero section of the cotangent bundle $T^*X$, and $\dim \Char(\cO_X)=n$.
\end{ejemplo}

The characteristic variety\index{Characteristic variety} $\Char(\cM) \subset T^*X$ of a coherent $\cD_X$--module is involutive (e.g. \cite{malgrange-involutivite-1978}, \cite{gabber-integrabilite}, \cite[App. B]{Granger-Maisonobe-cimpa-niza}). Hence, if $\cM$ is a non zero module, one has  $2n \geq \dim(\Char(\cM)) \geq n$. This is called Bernstein's inequality. It is proved in \cite{bernstein-ana-continuation-1972} for modules over the Weyl algebra.

\begin{definicion} A coherent $\cD_X$--module is said to be {\em holonomic}\index{Holonomic $\cD$--module} if either $\cM=\{0\}$ or $\dim (\Char(\cM)) = \dim\, X$.
\end{definicion}

The $\cD_X$--module $\cO_X$ is holonomic while $\cD_X$ is not holonomic as $\cD_X$--module. If $\cM=\frac{\cD_X}{\cI}$ where $\cI$ is a sheaf of non zero locally principal ideals (i.e. locally generated by a single differential operator) then $\cM$ is holonomic if and only if $\dim\, X=1$, since $\dim \Char(\frac{\cD_X}{\cI})=2n-1$. A central result in $\cD_X$--module theory is

\begin{teorema} {\rm \cite{Kashiwara-1977}}
If $D$ is a hypersurface in a complex manifold $X$ then the $\cD_X$--module $\cO_X(\star D)$ of meromorphic functions in $X$ with poles on $D$ is holonomic.
\end{teorema}

We finish this subsection with a result that we use later on (see Sect. \ref{deRhamcomplex-of-a-Dmodule}):

\begin{teorema}\label{teo:ext-de-un-mod-holonomo}
A non zero coherent $\cD_X$--module $\cM$ is holonomic if and only if $\fExt^i_{\cD_X}(\cM,\cD_X)=0$ for $i\neq n$. In this case, $\fExt^n_{\cD_X}(\cM,\cD_X)$ is a coherent holonomic right $\cD_X$--module.
\end{teorema} The proof follows from \cite[IV; Ths. 4.2.5, 4.2.6]{sem-gre-1975}, see also \cite[Ths. 7, 8]{Granger-Maisonobe-cimpa-niza}.

\subsection{The natural isomorphism $\Sym_{\cO_X}(\fDer(\cO_X)) \stackrel{\sim}{\rightarrow} \gr_{F^\bullet}(\DD_X)$}\label{symmetric-algebra} \label{sec:sim-der-gr-D}

There is a natural injective morphism of $\cO_X$--modules $$\fDer(\cO_X) \stackrel{\iota}{\rightarrow} \gr^{(1)}_{F^\bullet}(\DD) \subset  \gr_{F^\bullet}(\DD)$$ mapping a local section $\delta \in \fDer(\cO_X)$ to $\sigma(\delta)\in \gr^{(1)}_{F^\bullet}(\DD)$. Then, by the universal property of the symmetric algebra (see \cite[\S 6, Prop. 2]{Bourbaki-algebra-1-3}) there exists a unique morphism of $\cO_X$--algebras
\beq\label{kappa}
\kappa : \Sym_{\cO_X}(\fDer(\cO_X)) \rightarrow \gr_{F^\bullet}(\DD_X),
\eeq
extending the morphism $\iota: \fDer(\cO_X) \rightarrow  \gr_{F^\bullet}(\DD).$ Since $\gr_{F^\bullet} (\DD_X)$ is a graded $\cO_X$--algebra and $\iota(\fDer(\cO_X)) \subset \gr^{(1)}_{F^\bullet}(\DD)$, by [Rmk. 1, page 498, loc.cit.] the morphism $\kappa$ is graded so that we have, for any $k\geq 0$ a morphism of $\cO_X$--modules
$$\kappa_k : \Sym^{(k)}_{\cO_X}(\fDer(\cO_X)) \rightarrow \gr^{(k)}_{F^\bullet}(\DD_X).$$

By \cite[\S 6, Th. 1]{Bourbaki-algebra-1-3} and since $\fDer(\cO_X)$ is a locally free $\cO_X$--module of rank $n$, for any $x\in X$, $\Sym_{\cO_{X,x}}(\fDer(\cO_{X,x}))$
is canonically isomorphic to the polynomial algebra $\cO_{X,x}[T_1,\ldots,T_n]$, the canonical isomorphism being obtained by mapping $\partial_i$ to $T_i$ for $1\leq i \leq n$. Locally, for $\alpha = (\alpha_1,\ldots,\alpha_n)$ and $|\alpha| =\sum_i \alpha_i = k$, one has   $\kappa_k(\partial_1^{\alpha_1}\cdots \partial_n^{\alpha_n})= \sigma(\partial_1^{\alpha_1}\cdots \partial_n^{\alpha_n})$
where $\sigma(\,\,)$ means the principal symbol (see (\ref{basic-objects})). Then $\kappa_k$ is an isomorphism for any $k$ and $\kappa$ is a graded isomorphism. The isomorphism $\kappa$ is the intrinsic version of the isomorphism (\ref{grD-O-xi}).

\subsection{The Spencer complex\index{Spencer complex}} \label{Spencer-complex-subsection}

Let $X$ be a complex manifold of dimension $n\geq 1$.

\begin{definicion} \label{Spencer-complex} {\rm\cite[Chap I, (2.1)]{meb_formalisme}} The Spencer complex in $X$ is the following complex of left $\DD_X$--modules, denoted  $\Spb_X$ (or simply $\Spb$):
\begin{equation*}
0 \rightarrow \DD \otimes_{\cO_X} \bigwedge^n \fDer(\cO_X) \stackrel{\epsilon_{-n}}{\longrightarrow} \cdots \stackrel{\epsilon_{-2}}{\longrightarrow} \DD \otimes_{\cO_X} \fDer(\cO_X) \stackrel{\epsilon_{-1}}{\longrightarrow} \DD
\end{equation*} where the differential $\epsilon_\bullet$ is defined by
\begin{equation*}
\begin{split}
\epsilon_{-1}(P\otimes \delta) & = P\delta, \\
\epsilon_{-p}(P\otimes (\delta_1 \wedge \cdots \wedge \delta_p)) & = \sum_{i=1}^p (-1)^{i-1}P\delta_i \otimes (\delta_1 \wedge \cdots \widehat{\delta_i} \cdots \wedge \delta_p)) \, + \\
& \sum_{1\leq i < j \leq p}^{} (-1)^{i+j} P \otimes ([\delta_i,\delta_j] \wedge \delta_1 \wedge \cdots \widehat{\delta_i} \cdots  \widehat{\delta_j} \cdots \wedge \delta_p)\,\,\, (2\leq p\leq n).
\end{split}
\end{equation*}
We denote by $\widetilde{\Spb}_X$ (or simply $\widetilde{\Spb}$) the augmented complex $$\Sp_X^\bullet \stackrel{\epsilon_0}{\longrightarrow} \cO_X \rightarrow 0$$ where $\epsilon_0(P)=P(1)$.
\end{definicion}

\begin{nota}\label{SP-free-resolution-OX}
The complex $\Spb_X$ is a locally free resolution of the left $\cD_X$--module $\cO_X$. Since $\fDer(\cO_X)$ is a locally free $\cO_X$--module, each left $\cD_X$--module $$\Sp_X^p  := \DD_X \otimes_{\cO_X} \bigwedge^p \fDer(\cO_X)$$ is locally free (or rank $\binom{n}{p}$).

To prove that $\Spb_X$ is a resolution of $\cO_X$, or equivalently that $\widetilde{\Spb}_X$ is acyclic, we can proceed locally at each point $x\in X$. If we choose local coordinates $(x_1,\ldots,x_n)$ on $X$ around $x$, the partial derivatives $(\partial_1,\ldots,\partial_n)$ -- where $\partial_i = \frac{\partial}{\partial x_i}$ -- form a basis of the free $\cO_X$--module $\fDer(\cO_X)_x$. The differential $\epsilon_{-p,x}$ can be written as
$$ \epsilon_{-p,x}(P\otimes (\partial_{i_1} \wedge \cdots \wedge \partial_{i_p}))  = \sum_{j=1}^p (-1)^{j-1}P\partial_{i_j} \otimes (\partial_{i_1} \wedge \cdots \widehat{\partial_{i_j}} \cdots \wedge \partial_{i_p}).$$
So, the complex $\widetilde{\Spb}_{X,x}$ is a Koszul complex (which is denoted by $K(\partial_1, \ldots,\partial_n; \DD_X)$ in \cite[I.2]{Malgrange-cimpa-niza}. This complex is exact (see e.g. \cite{Frisch}). We give here a proof of the acyclicity of $\widetilde{\Spb}_{X,x}$ based on  \cite[Th. 3.1.2; Prop. 4.1.3]{calde_ens} (which proves a more general result and follows a suggestion of B. Malgrange \cite[I.2]{Malgrange-cimpa-niza}).

We consider a discrete increasing filtration $G^\bullet:=G^\bullet\left(\widetilde{\Spb}_{X}\right)$ on the complex $\widetilde{\Spb}_{X}$ (or more precisely on the complex $\widetilde{\Spb}_{X,x}$). The discrete filtration $G^\bullet$ is compatible with the differentials and the associated graded complex is exact. This implies that the complex $\widetilde{\Spb}_{X}$ is exact since the filtration $G^\bullet$ is discrete (i.e. $G^k=0$ if $k<0$).

The definition of $G^{\bullet}$ is as follows:  for $0 \leq p \leq n$ and $k\in \NN$, write $$G^{k,-p}:=G^k\left(\DD \otimes_{\cO_X} \bigwedge^p \fDer(\cO_X)\right) = F^{k-p}(\DD) \otimes _{\cO_X} \bigwedge^p \fDer(\cO_X)$$  where $F^\bullet(\DD)$ is the order filtration in $\DD$. We also write $G^{k,1}:=G^k(\cO_X)= \cO_X$ for all $k\in \NN$.

For each $k,p$ one has $\epsilon_{-p} (G^{k,-p}) \subset G^{k,-p+1}$ since for any $P\in F^{k-p}(\DD)$ one has $P\partial_i \in F^{k-p+1}(\DD)$. Thus, $G^{k,\bullet}$ is a complex of $\cO_X$--modules for any $k$ and $(G^{k,\bullet})_k$ is a discrete increasing exhaustive filtration (simply denoted by $G^\bullet$) of the complex $\widetilde{\Spb}_{X}$.

For each $0 \leq p \leq n$, the family $(G^{k,-p})_k$ is a discrete increasing exhaustive filtration of the $\DD$--module $\DD\otimes \bigwedge^p \fDer(\cO_X)$ whose associated graded $\cO_X$--module  $$\gr_{G^{\bullet,-p}}\left(\DD\otimes \bigwedge^p \fDer(\cO_X)\right)$$ is naturally isomorphic to
$$\gr_{F^\bullet}(\DD)[-p]\otimes \bigwedge^p \fDer(\cO_X).$$
The last isomorphism follows from the fact that  $\bigwedge^p \fDer(\cO_X)$ is a locally free $\cO_X$--module and therefore $\cO_X$--flat.

Then the associated graded complex $\gr_{G^\bullet}(\widetilde{\Spb})$ is

\begin{equation*}
0 \rightarrow \gr_{F^\bullet}(\DD)[-n] \otimes \bigwedge^n \fDer(\cO_X) \stackrel{\tau_{-n}}{\longrightarrow} \cdots \stackrel{\tau_{-1}}{\longrightarrow} \gr_{F^\bullet}(\DD)   \stackrel{\tau_{0}}\rightarrow \cO_X  \rightarrow 0
\end{equation*} where the differential $\tau_\bullet=\gr(\epsilon_\bullet)$ is acting by
\begin{equation*}
\begin{split}
\tau_0(Q)=Q_0,\\
\tau_{-1}(Q\otimes \partial_i) & = Q\sigma(\partial_i), \\
\tau_{-p}(Q\otimes (\partial_{i_1} \wedge \cdots \wedge \partial_{i_p})) & = \sum_{j=1}^p (-1)^{j-1}Q\sigma(\partial_{i_j}) \otimes (\partial_{i_1} \wedge \cdots \widehat{\partial_{i_j}} \cdots \wedge \partial_{i_p})) \,\,\, (2\leq p\leq n),
\end{split}
\end{equation*}
where the tensor product is taken over $\cO_X$,  $Q_0$ is the $0-th$ homogenous component of $Q\in \gr_{F^\bullet}(\DD)$ and $\sigma(\,)$ is the principal symbol of the corresponding differential operator (see (\ref{basic-objects})).
So, the last complex is nothing but the augmented Koszul complex with respect to the regular sequence $\{\sigma(\partial_1), \ldots, \sigma(\partial_n)\}$ in the commutative ring $\gr_{F^\bullet}(\DD)$. Hence this complex is acyclic.
\end{nota}

\subsection{The de Rham complex of a $\cD_X$--module\index{de Rham complex of a $\cD$--module}}\label{deRhamcomplex-of-a-Dmodule}

We follow here \cite[Chap. I (2.6)]{meb_formalisme}. Let $\MM$ be a left $\DD_X$--module. The left exact functor  $$\fHom_{\DD_X}(\cO_X, *): \Mod(\DD_X) \longrightarrow \Mod(\CC_X)$$ can be derived to give a functor
$$\RR \fHom_{\DD_X}(\cO_X, *) : D^b(\DD_X) \longrightarrow D^b(\CC_X).$$ Here $\Mod(A_X)$ stands for the category of (left) $A_X$--modules and $D^b(A_X)$ for the derived category of complexes of $A_X$--modules with bounded cohomology.

\begin{definicion}\label{de-Rham}
Let $\MM$ be a left $\DD_X$--module. The complex $\RR \fHom_{\DD_X}(\cO_X, \MM)$ in $D^b(\CC_X)$ is called the de Rham complex of $\MM$. It is denoted by $\DR(\MM)$.
\end{definicion}

Let $\MM$ be a left $\DD_X$--module. Since $\MM$ carries an integrable connection, there is a natural morphism of sheaves of $\CC_X$--vector spaces $$\nabla : \MM \longrightarrow \Omega^1_X \otimes_{\cO_X} \MM $$
given locally by $$\nabla(m) = \sum_{i=1}^n  dx_i \otimes \partial_i(m).$$

\begin{proposicion} \label{Omega-DR}
{\rm \cite[Lemme (2.6.3)]{meb_formalisme}} For any left $\DD_X$--module $\MM$, the complex $\DR(\MM)$ can be represented by the complex
$$ \Omega^\bullet_X(\MM):=0 \rightarrow \MM \stackrel{\nabla}{\longrightarrow} \Omega^1_X \otimes _{\cO_X} \MM  \stackrel{\nabla}{\longrightarrow}  \cdots \stackrel{\nabla}{\longrightarrow} \Omega^n_X \otimes _{\cO_X} \MM  \rightarrow 0$$
concentrated in degrees $[0,n]$, where $\nabla(\omega \otimes m) = \omega\wedge  \nabla(m)  + (-1)^{\deg(\omega)} d\omega \otimes  m.$
\end{proposicion}

\begin{prueba} By Remark \ref{SP-free-resolution-OX} the complex $\widetilde{\Sp}^\bullet$ is a locally free resolution of the left $\DD_X$--module $\cO_X$. Thus $\DR(\MM)$ is represented by $\fHom_{\cD_X} (\Sp^\bullet,\MM)$. Moreover, for $p=0,\ldots,n$, the $\cO_X$--modules $\fHom_{\cD_X}(\cD_X\otimes_{\cO_X} \bigwedge^p \fDer(\cO_X), \MM)$ and $\Omega^p_X \otimes \cM$ are isomorphic. And there exists a natural quasi-isomorphism, in $D^b(\CC_X)$, from $\fHom_{\cD_X}(\cD_X\otimes_{\cO_X} \bigwedge^\bullet  \fDer(\cO_X), \MM)$ to $\Omega^\bullet_X \otimes \MM$.
\end{prueba}

Following \cite[Ch. 1, (4.1)]{meb_formalisme}, if $\cM$ is a complex in the derived category $D^b_c(\cD_X)$ (complexes of left $\cD_X$--modules, with bounded and coherent cohomology), the complex $$\RR \fHom_{\DD_X}(\cM,\cD_X)$$ of right $\cD_X$--modules has bounded and coherent cohomology.

\begin{definicion}\label{Dual-D-modules}
Let $\cM$ be a complex in $D^b_c(\cD_X)$. The dual\index{$\cD$--module dual} of $\cM$ is the complex $\cM^*$ in $D^b(\cD_X)_c$ defined by $$\cM^* = \fHom_{\cO_X}(\Omega^n_X, \RR \fHom_{\DD_X}(\cM,\cD_X))[n]$$ where $n=\dim X$.
\end{definicion}

The dual $\cM^*$ is also denoted $\dsD(\cM)$. The dual $\cM^*$ of a holonomic $\cD_X$--module $\cM$ is also holonomic (see \cite[Ch. I, (4.1)]{meb_formalisme} and Th. \ref{teo:ext-de-un-mod-holonomo}). The dual $\cO_X^*$ is naturally isomorphic to $\cO_X$ and there is a natural isomorphism in the derived category $D^b(\CC_X)$ $$\DR(\cM) \stackrel{\sim}{\rightarrow} \RR \fHom_{\cD_X}(\cM^*,\cO_X).$$ The last complex is called the (holomorphic)  solution complex of $\cM^*$.

\subsection{Grothendieck's Comparison Theorem revisited}\label{Gro-revisited}

The statement of (a version of) Grothendieck's Comparison Theorem is as follows (see Sect. \ref{Grothendieck-comparison-theorem}):

\begin{teorema}\index{Grothendieck's Comparison Theorem} {\rm (\cite[Th. 2]{GRO})} \label{teo:Gro-comp-th-revisited}
If $D$ is a divisor in the complex manifold $X$, and $j : U:=X\ssm D \hookrightarrow X$ is the inclusion, then the de Rham morphism
\beq\label{de-Rham-morphism-Gro-revisited} \Omega^\bullet_X(\star D) \longrightarrow {\RR} j_*\CC_U
\eeq
is a quasi-isomorphism. \end{teorema}

\begin{nota}
If $D$ is a normal crossing divisor, the result is due to Atiyah-Hodge \cite{AH55}. The proof of Grothendieck uses Hironaka's resolution of singularities to reduce the general case to the normal crossing divisor one.

Theorem 2 in \cite{GRO} is more general. It holds for a reduced complex analytic space $X$, an analytic closed subset  $D\subset X$,  assuming $U:=X\ssm D$ is non singular and dense in $X$, and that $U$ can be defined locally by one
equation.
\end{nota}

Mebkhout (\cite[Chap. 2, \S 2]{meb_formalisme})\footnote{Following Mebkhout's Ph.D. Thesis} interprets Grothendieck's comparison theorem as the {\em regularity} of the $\cD_X$--module $\cO_X$. This regularity is equivalent to the fact that, in the derived category $D^b(\CC_X)$,
the natural morphism $$\DR(\cO_X(\star D)) \rightarrow \DR(\RR j_*j^{-1} \cO_X)$$ is a quasi-isomorphism for any divisor $D$. This last morphism is the de Rham morphism (\ref{de-Rham-morphism-Gro-revisited}) since $\DR(\cO_X(\star D))$ equals $\Omega_X^\bullet(\star D)$ and $\DR(\RR j_*j^{-1} \cO_X)$ is quasi-isomorphic to $\RR j_*j^{-1} \DR(\cO_X) = \RR j_*j^{-1}\Omega^\bullet_X \simeq \RR j_* \CC_U$, by Poincar\'e Lemma.

\section{Free divisors and logarithmic $D$--modules}\label{sec:free-divisors-log-D-modules}

\subsection{The sheaf of logarithmic differential operators with respect to a free divisor\index{Logarithmic differential operators}}\label{sub:sheaf-V}

Let $X$ be a complex analytic manifold of dimension $n$ and $D\subset X$ a hypersurface.
Let us denote $U= X \ssm D$ and $j:U\hookrightarrow X$ the corresponding open immersion.

When $D$ is a free divisor, there is a nice sheaf of subrings of the sheaf $\DD_X$ of linear differential operators in $X$, namely the sheaf of {\em logarithmic differential operators} with respect to $D$, denoted by $\VV^D_X$. It is the
sheaf of subrings of $\DD_X$ with stalks
$$
\VV^D_{X,x} = \{P\in \DD_{X,x}\ |\ P(\cJ_x^j) \subset \cJ_x^j\ \forall j\geq 0\}.
$$
Moreover, $\VV_X$ is a sheaf of filtered rings, with the induced filtration by the order filtration in $\DD_X$, whose graded ring is commutative and the canonical map
\begin{equation} \label{eq:kappa-log}
\Sym_{\OO_X} \fDer(-\log D) \longrightarrow \gr \VV^D_X
\end{equation}
is an isomorphism of commutative graded $\OO_X$-algebras (see \cite[Cor. 2.1.6]{calde_ens}; compare with Sect. \ref{sec:sim-der-gr-D}). As a consequence, $\VV^D_X$
is generated by $\OO_X$ and $\fDer(-\log D)$, and $\VV^D_X$ is the enveloping algebra of $\fDer(-\log D)$ considered as a Lie algebroid (cf. \cite[\S (2.1)]{chem_99}). From there we deduce that $\VV^D_X$ is a left and right coherent sheaf of rings (one can proceed as in the case of $\DD_X$, Sect. \ref{basic-objects}, or as in \cite[Th. 1.2.5]{Bjork_AD}) and its stalk at each point of $X$ is a left and right Noetherian ring of finite global homological dimension $\leq 2n$ (cf. \cite[App. IV, Prop. 4.14 and Th. 5.1]{Bjork_AD}).

\subsection{The logarithmic Spencer complex\index{Logarithmic Spencer complex}}\label{log-Spencer-complex}

From now on, we will assume that $D\subset X$ is a free divisor in a complex manifold $X$ of dimension $n\geq 1$.

\begin{definicion}  \label{def-log-Spencer-complex}{\rm\cite[Def. 3.1.1]{calde_ens}} The logarithmic Spencer complex associated with $D\subset X$ is the following complex of left $\VV_X^D$--modules, denoted  $\Sp^\bullet(\log D)$:
\begin{equation*}
0 \rightarrow \VV_X^D \otimes_{\cO_X} \bigwedge^n \fDer(-\log D) \stackrel{\epsilon_{-n}}{\longrightarrow} \cdots \stackrel{\epsilon_{-2}}{\longrightarrow} \VV_X^D \otimes_{\cO_X} \fDer(-\log D) \stackrel{\epsilon_{-1}}{\longrightarrow} \VV_X^D
\end{equation*} where the differential $\epsilon_\bullet$ is defined by
\begin{equation*}
\begin{split}
\epsilon_{-1}(P\otimes \delta) & = P\delta, \\
\epsilon_{-p}(P\otimes (\delta_1 \wedge \cdots \wedge \delta_p)) & = \sum_{i=1}^p (-1)^{i-1}P\delta_i \otimes (\delta_1 \wedge \cdots \widehat{\delta_i} \cdots \wedge \delta_p)) \, + \\
& \sum_{1\leq i < j \leq p}^{} (-1)^{i+j} P \otimes ([\delta_i,\delta_j] \wedge \delta_1 \wedge \cdots \widehat{\delta_i} \cdots  \widehat{\delta_j} \cdots \wedge \delta_p)\,\,\, (2\leq p\leq n).
\end{split}
\end{equation*}
We denote by $\widetilde{\Sp}^\bullet(\log D)$ the augmented complex $$\Sp^\bullet(\log D) \stackrel{\epsilon_0}{\longrightarrow} \cO_X$$ where $\epsilon_0(P)=P(1)$.
\end{definicion}


For a free divisor $D$ the complex $\Sp^\bullet(\log D)$ is a locally free resolution of the left $\VV_X^D$-module $\OO_X$ \cite[Th. 3.1.2]{calde_ens} (this is a particular case of Proposition \ref{log-Spencer-resolution-E}).

\begin{nota} The logarithmic Spencer complex should be compared with the Spencer complex (see Definition \ref{Spencer-complex}).

The logarithmic Spencer complex generalizes the one given in \cite[App. A (A.4)]{esnault-vieweg-1986} for a normal crossing divisor $D\subset X$. Actually, this definition is a sheaf version of the Rinehart complex of a Lie-Rinehart algebra (see \cite{rine-63}).
\end{nota}

A {\em logarithmic connection}\index{Logarithmic connection} (with respect to $D$) on an $\OO_X$-module $\calE$ is a $\CC$-linear map $$\nabla: \calE \tos \Omega^1_X(\log D)\otimes_{\OO_X}\calE  $$ satisfying the Leibniz rule, i.e. $\nabla (a e) = da \otimes e + a \nabla(e) $ for any holomorphic function $a$ and any local section $e$ of $\calE$.

As in the classical case (see Sect. \ref{deRhamcomplex-of-a-Dmodule}), such a $\nabla$ may be extended to a family of $\CC$-linear maps $$\nabla^p: \Omega_X^p(\log D)\otimes_{\OO_X}\calE \tos \Omega_X^{p+1}(\log D)\otimes_{\OO_X}\calE,$$ with $\nabla^p(\alpha \otimes e) = d \alpha \otimes e + (-1)^p \alpha \wedge \nabla(e)$ for any logarithmic $p$-form $\alpha$ and any local section $e$ of $\calE$.

We say that the logarithmic connection $\nabla$ is {\em integrable}\index{Integrable logarithmic connection} if $\nabla^{p+1} \circ \nabla^p=0$ for all $p\geq 0$. In such a case, the {\em logarithmic de Rham complex} of $(\calE,\nabla)$ is by definition the complex (of sheaves of complex vector spaces)
$$ \Omega^\bullet_X(\log D)(\calE) :=
\calE \tos \Omega^1_X(\log D)\otimes_{\OO_X}\calE \tos \Omega^2_X(\log D)\otimes_{\OO_X}\calE \tos \cdots \tos  \Omega^n_X(\log D)\otimes_{\OO_X}\calE,
$$
where $\calE$ is placed in degree $0$.

Any logarithmic connection $\nabla$ on an $\OO_X$-module $\calE$ gives rise to an action of logarithmic vector fields $\fDer(-\log D)$ on $\calE$
$$ (\delta,e) \longmapsto \nabla_\delta(e) = \langle \delta,\nabla(e)\rangle$$
for any logarithmic vector field $\delta$ and any local section $e$ of $\calE$,
where $\langle \delta,\nabla(e)\rangle$ is induced by the contraction of logarithmic 1-forms by logarithmic vector fields.

Obviously, the exterior derivative $d:\OO_X \tos \Omega^1_X(\log D)$ is a logarithmic connection which is integrable, and $\Omega^\bullet_X(\log D)(\OO_X) = \Omega^\bullet_X(\log D)$.

When $D$ is free, the following result holds,  by essentially the same proof as in the classical case.

\begin{proposicion} Assume that $D$ is a free divisor and
let $\nabla$ be a logarithmic connection on an $\OO_X$-module $\calE$. The following properties are equivalent:
\begin{enumerate}
\item[(i)] $\nabla$ is integrable.
\item[(ii)] $\nabla_{[\delta,\delta']} = [\nabla_\delta,\nabla_{\delta'}]$ for all logarithmic vector fields $\delta,\delta'$.
\end{enumerate}
\end{proposicion}

Since the sheaf of logarithmic differential operators $\VV_X^D$ is the enveloping algebra of the Lie algebroid $\fDer(-\log D)$ provided that $D$ is a free divisor, we obtain the following corollary.

\begin{corolario} Assume that $D$ is a free divisor and let $\calE$ be an $\OO_X$-module. The following data are equivalent:
\begin{enumerate}
\item[(a)] An integrable logarithmic connection $\nabla$ on $\calE$.
\item[(b)] A structure of left $\VV_X^D$-module on $\calE$ extending its $\OO_X$-module structure.
\end{enumerate}
Moreover, the action of a logarithmic derivation $\delta$ on $\calE$ in (b) is given by $\delta\cdot e = \nabla_\delta(e)$ for each local section $e$ of $\calE$.
\end{corolario}

From now on, we assume that $D$ is a free divisor.

Any locally free $\OO_X$-module of finite rank endowed with an integrable logarithmic connection will be called an ILC, for short.
Examples of ILCs are the invertible $\OO_X$-modules $\OO_X(k D)$, $k\in \ZZ$, which carry an obvious structure of left $\VV_X^D$-module.

We define the logarithmic Spencer complex of an arbitrary left $\VV_X^D$-module in the following way.

\begin{definicion}  \label{log-Spencer-complex-E}{\rm\cite{calde_nar_fourier}} Let $\calE$ be a left $\VV_X^D$-module.
The logarithmic Spencer complex\index{Logarithmic Spencer complex} of $\calE$ (with respect to $D\subset X$) is the following complex of left $\VV_X^D$--modules, denoted  $\Sp^\bullet(\log D)(\calE)$:
\begin{equation*}
0 \rightarrow \VV_X^D \otimes_{\cO_X} \bigwedge^n \fDer(-\log D) \otimes_{\cO_X} \calE \stackrel{\epsilon_{-n}}{\longrightarrow} \cdots \stackrel{\epsilon_{-2}}{\longrightarrow} \VV_X^D \otimes_{\cO_X} \fDer(-\log D)\otimes_{\cO_X} \calE \stackrel{\epsilon_{-1}}{\longrightarrow} \VV_X^D\otimes_{\cO_X} \calE
\end{equation*} where the differential $\epsilon_\bullet$ is defined by
\begin{equation*}
\begin{split}
\epsilon_{-1}(P\otimes \delta \otimes e)  = P\delta \otimes e - P\otimes \delta e,& \\
\epsilon_{-p}(P\otimes (\delta_1 \wedge \cdots \wedge \delta_p)) \otimes e & =\\
\sum_{i=1}^p (-1)^{i-1}P\delta_i \otimes (\delta_1 \wedge \cdots \widehat{\delta_i} \cdots \wedge \delta_p))\otimes e \, - \,
\sum_{i=1}^p (-1)^{i-1}P \otimes (\delta_1 \wedge \cdots \widehat{\delta_i} \cdots \wedge \delta_p))\otimes \delta_i e \, +&
\\
 \sum_{1\leq i < j \leq p}^{} (-1)^{i+j} P \otimes ([\delta_i,\delta_j] \wedge \delta_1 \wedge \cdots \widehat{\delta_i} \cdots  \widehat{\delta_j} \cdots \wedge \delta_p) \otimes e\,\,\, (2\leq p\leq n).&
\end{split}
\end{equation*}
We denote by $\widetilde{\Sp}^\bullet(\log D)(\calE)$ the augmented complex $$\Sp^\bullet(\log D)(\calE) \stackrel{\epsilon_0}{\longrightarrow} \calE$$ where $\epsilon_0(P\otimes e)=P e$.
\end{definicion}

Note that for $\calE=\OO_X$ we have $\Sp^\bullet(\log D)(\OO_X) = \Sp^\bullet(\log D)$.

\begin{proposicion}\label{log-Spencer-resolution-E} For any ILC $\calE$, the logarithmic Spencer complex of $\calE$ (with respect to $D\subset X$) $\Sp^\bullet(\log D)(\calE)$ is a locally free resolution of $\calE$.
\end{proposicion}

The proof of this proposition is similar to the proof in Remark \ref{SP-free-resolution-OX}. Namely, we consider the discrete increasing filtration
$G^\bullet:=G^\bullet\left(\widetilde{\Sp}^\bullet(\log D)(\calE)\right)$ on the complex $\widetilde{\Sp}^\bullet(\log D)(\calE)$ given, for $0 \leq p \leq n$ and $k\in \NN$, by
$$G^{k,-p}:=G^k\left(\VV_X^D  \otimes_{\cO_X} \bigwedge^p \fDer(-\log D) \otimes_{\cO_X} \calE\right) = F^{k-p}(\VV_X^D ) \otimes _{\cO_X} \bigwedge^p \fDer(-\log D) \otimes_{\cO_X} \calE$$  where $F^\bullet(\VV_X^D )$ is the filtration induced by the order filtration in $\DD$. We also write $G^{k,1}:=G^k(\calE)= \calE$ for all $k\in \NN$.
The associated graded complex turns out to be (locally) the tensor product over $\cO_X$ of
$\calE$ and the augmented Koszul complex with respect to
$\{\sigma(\delta_1),\dots,\sigma(\delta_n)\} \subset \gr \VV_X^D$, $\{\delta_1,\dots,\delta_n\}$ being a local $\cO_X$-basis of $\fDer(-\log D)$,
which is exact since $\calE$ is locally free over $\cO_X$ and $\gr \VV_X^D$ is (locally) a polynomial ring over $\cO_X$ in the variables $\{\sigma(\delta_1),\dots,\sigma(\delta_n)\}$ (see the isomorphism (\ref{eq:kappa-log})).
\medskip

By using the logarithmic Spencer resolution \ref{log-Spencer-resolution-E}, as in the case of $\DD_X$ (see Sect. \ref{Spencer-complex-subsection}), we obtain a canonical isomorphism of complexes of sheaves of complex vector spaces
$$
\fHom_{\VV_X^D}(\SP^\bullet(\log D),\calE) \simeq \Omega^\bullet_X(\log D)(\calE),
$$
for any left $\VV_X^D$-module $\calE$,
and so an isomorphism in the derived category $\RR\fHom_{\VV_X^D}(\OO_X,\calE) \simeq \Omega^\bullet_X(\log D)(\calE)$. This is completely similar to Proposition \ref{Omega-DR}.

\subsection{A $\DD$--module criterion for LCT\index{$\DD$--module criterion for LCT}}\label{D-module-criterion-for-LCT}

Grothendieck's comparison theorem (see (\ref{Grothendieck-comparison-theorem}) and (\ref{Gro-revisited})) tells us that the natural map $\Omega^\bullet_X(\star D) \tos \RR j_* \CC_{X\ssm D} $ obtained by composition of the adjunction map $\Omega^\bullet_X(\star D) \tos \RR j_* j^{-1} \Omega^\bullet_X(\star D) = \RR j_*\Omega^\bullet_{X\ssm D} = j_*\Omega^\bullet_{X\ssm D}$ with the inverse of the induced map by the Poincar\'e quasi-isomorphism $\CC_{X\ssm D} \tos \Omega^\bullet_{X\ssm D}$ is an isomorphism in the derived category of sheaves of complex vector spaces.
Definition \ref{def:LCT-holds-for-D} tells us that LCT holds for $D$ if the natural morphism
\begin{equation} \label{eq:log-compar}
\Omega^\bullet_X(\log D) \hookrightarrow \Omega^\bullet_X(\star D)
\end{equation}
is a quasi-isomorhism. Let us explain how this map
can be interpreted in terms of $\DD$--module theory in the case of free divisors.

On one hand, we have canonical isomorphisms (in the derived category)
$$
\Omega^\bullet_X(\log D) \simeq \RR\fHom_{\VV_X^D}(\OO_X,\OO_X) \simeq \RR\fHom_{\DD_X}\left(\DD_X \Lotimes_{\VV_X^D} \OO_X,\OO_X\right).
$$
By taking $\DD_X$-duals\index{$\cD$--module dual} $\dsD(-)$ (see Def. \ref{Dual-D-modules}) we obtain
\begin{eqnarray*}
&\displaystyle
\RR\fHom_{\DD_X}\left(\DD_X \Lotimes_{\VV_X^D} \OO_X,\OO_X\right) \simeq \RR\fHom_{\DD_X}\left(\dsD\, \OO_X, \dsD \left(\DD_X \Lotimes_{\VV_X^D} \OO_X\right)\right) \simeq &\\
& \displaystyle \RR\fHom_{\DD_X}\left( \OO_X, \dsD \left(\DD_X \Lotimes_{\VV_X^D} \OO_X\right)\right) = \DR \left( \dsD \left(\DD_X \Lotimes_{\VV_X^D} \OO_X\right)\right),
\end{eqnarray*}
and by using the duality formula in \cite[Cor. 3.1.2]{calde_nar_fourier} (see also
\cite[Th. (A.32)]{nar_symmetry_BS}), we obtain
$$
\dsD \left(\DD_X \Lotimes_{\VV_X^D} \OO_X\right) \simeq \DD_X \Lotimes_{\VV_X^D} \OO_X(D),
$$
and so
\begin{equation} \label{eq:2}
\Omega^\bullet_X(\log D) \simeq \cdots \simeq \DR\left(\DD_X \Lotimes_{\VV_X^D} \OO_X(D)\right).
\end{equation}
On the other hand we have another canonical isomorphism (see Proposition \ref{Omega-DR})
\begin{equation} \label{eq:3}
\Omega^\bullet_X(\star D) = \Omega^\bullet_X(\cO_X(\star D)) \simeq \RR\fHom_{\DD_X}(\OO_X,\OO_X(\star D)) = \DR(\OO_X(\star D)).
\end{equation}

The point is that the inclusion $\OO_X(D) \subset \OO_X(\star D)$ induces a canonical left $\DD_X$-linear map
\begin{equation} \label{eq:4}
\varrho:\DD_X \Lotimes_{\VV_X^D} \OO_X(D) \longrightarrow \OO_X(\star D),
\end{equation}
and we have the following theorem (\cite[Cor. 4.2]{calde_nar_fourier} and \cite[Cor. 1.7.2]{nar_comp_08}):

\begin{teorema} \label{th:D-crit-LCT}
 Under the above hypotheses, the following properties hold:
\begin{enumerate}
\item[(1)] The logarithmic comparison map (\ref{eq:log-compar}) corresponds to the map (\ref{eq:4}) under the $\DR$ functor and the canonical isomorphisms (\ref{eq:2}) and (\ref{eq:3}).
\item[(2)] The following properties are equivalent:
\begin{enumerate}
\item[(2-1)] The logarithmic comparison theorem holds for $D$, i.e. the map (\ref{eq:log-compar}) is a quasi-isomorphism.
\item[(2-2)] The map (\ref{eq:4}) is an isomorphism (in the derived category of $\DD_X$-modules).
\item[(2-3)] The canonical map
$$j_! \CC_{X\ssm D} \longrightarrow \Omega^\bullet_X(\log D)(\cO_X(-D))$$
is a quasi-isomorphism.
\end{enumerate}
\end{enumerate}
\end{teorema}

\begin{nota} \label{rem:comple_D-crit-LCT} (a) Let us notice that property (2-3) above comes from taking Grothendieck-Verdier duals. Namely, $\left(\RR j_* \CC_{X\ssm D}\right)^\vee \simeq j_! \CC_{X\ssm D}$ and
\begin{eqnarray*}
&\displaystyle
\left( \Omega^\bullet_X(\log D)\right)^\vee \simeq
\left(\DR\left(\DD_X \Lotimes_{\VV_X^D} \OO_X(D)\right)\right)^\vee \stackrel{(\star)}{\simeq}
\RR\fHom_{\DD_X}\left(\DD_X \Lotimes_{\VV_X^D} \OO_X(D),\OO_X\right) \simeq &
\\
&\displaystyle
\RR\fHom_{\VV_X^D}\left( \OO_X(D),\OO_X\right) \simeq \RR\fHom_{\VV_X^D}\left(\OO_X^*, \OO_X(D)^*\right) \simeq
&
\\
&\displaystyle
\RR\fHom_{\VV_X^D}\left(\OO_X, \OO_X(-D)\right)\simeq \Omega^\bullet_X(\log D)(\cO_X(-D)),
\end{eqnarray*}
where we have used Mebkhout local duality formula in $(\star)$ (see \cite[Chap. I, (4.3)]{meb_formalisme}; see also \cite{narLDT}).

\noindent (b) Let us also notice that the complex $\Omega^\bullet_X(\log D)(\cO_X(-D))$ is a subcomplex of $\Omega^\bullet_X$, since locally $\Omega^\bullet_X(\log D)(\cO_X(-D)) = f\, \Omega^\bullet_X(\log D)$ with $f=0$ a reduced local equation of $D$, and $f\, \Omega^p_X(\log D) \subset \Omega^p_X$ for $p=0,\dots,n$. Moreover,
\begin{eqnarray*}
&\displaystyle
\Omega^\bullet_X(\log D)(\cO_X(-D)) \simeq \cdots \simeq \RR\fHom_{\DD_X}\left(\DD_X \Lotimes_{\VV_X^D} \OO_X(D),\OO_X\right) \simeq
&
\\
&\displaystyle
 \RR\fHom_{\DD_X}\left(\dsD\, \OO_X,\dsD\left(\DD_X \Lotimes_{\VV_X^D} \OO_X(D) \right)\right) \simeq
 \RR\fHom_{\DD_X}\left(\OO_X,\DD_X \Lotimes_{\VV_X^D} \OO_X \right) = \DR \left(\DD_X \Lotimes_{\VV_X^D} \OO_X \right)
\end{eqnarray*}
and we can prove that the inclusion
$\Omega^\bullet_X(\log D)(\cO_X(-D)) \hookrightarrow \Omega^\bullet_X$ comes from the map of $\DD_X$-modules
\beq \label{eq:D-otimes-V-O-O}
\DD_X \Lotimes_{\VV_X^D} \OO_X  \longrightarrow \OO_X,\quad P \otimes a \longmapsto P(a),
\eeq
by applying the $\DR(-)$ functor.
\end{nota}

\begin{nota}\label{lctconcrete}Let us understand in concrete terms the significance of the map (\ref{eq:4}) being an isomorphism.
First of all, (\ref{eq:4}) is an isomorphism if and only if it is so stalkwise, and obviously $\varrho_p$ is an isomorphism for each $p\in X \ssm D$. Let us take a point $p\in D$, a reduced local equation $f\in \OO_{X,p}$ of $D$ and a basis $\delta_1,\dots,\delta_n\in \fDer(-\log D)_p$ with $\delta_i(f) = \alpha_i f$, $\alpha_i\in \OO_{X,p}$, $i=1,\dots,n$.
On the other hand, $\OO_{X,p}(D)$ is generated as $\VV^D_{X,p}$-module by $f^{-1}$ and the kernel of
$$
P \in \VV^D_{X,p} \longmapsto P(f^{-1}) \in \OO_{X,p}(D)
$$
is the left ideal generated by $\delta_1 + \alpha_1,\dots,\delta_n+\alpha_n$, i.e.
$$\OO_{X,p}(D) \simeq \VV^D_{X,p}/\calV^D_{X,p} \langle \delta_1 + \alpha_1,\dots,\delta_n+\alpha_n\rangle.$$
Consequently the map $\varrho_p:\DD_{X,p} \Lotimes_{\VV^D_{X,p}} \OO_{X,p}(D) \longrightarrow \OO_{X,p}(\star D)$ is an isomorphism (in the derived category of $\DD_{X,p}$-modules) if and only if the following properties hold:
\begin{enumerate}
\item[(i)] The complex $\DD_{X,p} \Lotimes_{\VV^D_{X,p}} \OO_{X,p}(D)$ is exact in cohomological degrees $\neq 0$, and
\item[(ii)] $\OO_{X,p}(\star D) \simeq \DD_{X,p}/\DD_{X,p} \langle \delta_1 + \alpha_1,\dots,\delta_n+\alpha_n\rangle$.
\end{enumerate}
The isomorphism in (ii) comes from the map
$$
P \in \DD_{X,p} \longmapsto P(f^{-1}) \in \OO_{X,p}(\star D),
$$
and so property (ii) is equivalent to
\begin{enumerate}
\item[(ii-1)] $\OO_{X,p}(\star D)$ is generated as $\DD_{X,p}$-module by $f^{-1}$, and
\item[(ii-2)] The $\DD_{X,p}$-annihilator of $f^{-1} \in \OO_{X,p}(\star D)$ is the left ideal generated by $\delta_1 + \alpha_1,\dots,\delta_n+\alpha_n$.
\end{enumerate}
Notice that property (ii-2) is equivalent to the fact that the $\DD_{X,p}$-annihilator of $f^{-1}$ is generated by order $1$ operators. From \cite{torrelli_2004} we know that the last property implies that the $b$-function of the germ $f$ has no integer roots strictly less than $-1$, and so property (ii-2) implies property (ii-1) (by the Bernstein functional equation). We conclude that the map
$\varrho_p$ is an isomorphism if and only if the following properties hold:
\begin{enumerate}
\item[(i)] The complex $\DD_{X,p} \Lotimes_{\VV^D_{X,p}} \OO_{X,p}(D)$ is exact in cohomological degrees $\neq 0$, and
\item[(ii-2)] The $\DD_{X,p}$-annihilator of $f^{-1} \in \OO_{X,p}(\star D)$ is the left ideal generated by $\delta_1 + \alpha_1,\dots,\delta_n+\alpha_n$.
\end{enumerate}
\end{nota}

Let us recall the following definition \cite{cas_ucha_stek}.

\begin{definicion} \label{def:Spencer}
A free divisor $D\subset X$ is called {\em Spencer}\index{Spencer free divisor} if the complex $\DD_X \Lotimes_{\VV_X^D} \OO_X$ is concentrated in cohomological degree $0$ and $\DD_X \otimes_{\VV_X^D} \OO_X$ is holonomic.
\end{definicion}

If the the LCT holds for $D$, then $\dsD (\OO_X(\star D)) \simeq \dsD \left( \DD_X \Lotimes_{\VV_X^D} \OO_X(D) \right) \simeq \DD_X \Lotimes_{\VV_X^D} \OO_X$, and so $D$ is Spencer and $\dsD (\OO_X(\star D)) \simeq  \DD_X \otimes_{\VV_X^D} \OO_X$ (since both $\OO_X(\star D)$ and its dual\index{$\cD$--module dual} are holonomic).
\medskip

To check whether the complex $\DD_X \Lotimes_{\VV_X^D} \OO_X$ is concentrated in cohomological degree $0$, we can use the logarithmic Spencer complex  $\Sp^\bullet(\log D)$, which is a locally free resolution of $\cO_X$ as a left $\VV_X^D$-module. So we have to check whether the complex $\DD_X \otimes_{\VV_X^D} \Sp^\bullet(\log D)$ is concentrated in cohomological degree $0$ or not.

For this, we can consider again the discrete increasing filtration
$G^\bullet:=G^\bullet\left(\DD_X \otimes_{\VV_X^D}\Sp^\bullet(\log D)\right)$ on the complex $\DD_X \otimes_{\VV_X^D}\Sp^\bullet(\log D)$ given by, for $0 \leq p \leq n$ and $k\in \NN$,
$$G^{k,-p}:=G^k\left(\DD_X   \otimes_{\cO_X} \bigwedge^p \fDer(-\log D) \right) = F^{k-p}(\DD_X ) \otimes _{\cO_X} \bigwedge^p \fDer(-\log D) $$  where $F^\bullet(\DD_X )$ is the filtration by the order filtration in $\DD$.
The associated graded complex turns out to be (locally) the Koszul complex with respect to
$\{\sigma(\delta_1),\dots,\sigma(\delta_n)\} \subset \gr \DD_X$, $\{\delta_1,\dots,\delta_n\}$ being a local $\cO_X$-basis of $\fDer(-\log D)$, but in general, this is not a regular sequence\footnote{Remember that $\{\sigma(\delta_1),\dots,\sigma(\delta_n)\}$ is a regular sequence in $\gr \VV_X^D$ since this ring is (locally) a polynomial ring in the variables $\{\sigma(\delta_1),\dots,\sigma(\delta_n)\}$ with coefficients in $\cO_X$.} and so $\gr_{G^\bullet} \left(\DD_X \otimes_{\VV_X^D}\Sp^\bullet(\log D)\right)$ is not concentrated in cohomological degree $0$. This fact motivates the following definition \cite[Def. 4.1.1]{calde_ens}.

\begin{definicion} \label{def:Koszul}
We say that a free divisor $D\subset X$ is {\em Koszul}\index{Koszul free divisor} at a point $p\in D$ if for some (and hence for any) local basis $\{\delta_1,\dots,\delta_n\}$ of $\fDer(-\log D)_p$, the symbols $\{\sigma(\delta_1),\dots,\sigma(\delta_n)\}$ form a regular sequence in $\gr \DD_{X,p}$; and we say that $D$ is {\em Koszul} if it so at any point $p\in D$.
\end{definicion}

Let us notice that the Koszul property for a free divisor $D\subset X$ is equivalent to saying that the complex  $\gr \DD_X \Lotimes_{\gr\VV_X^D} \OO_X$ is concentrated in cohomological degree $0$. In that case, the module $\gr \DD_X \otimes_{\gr \VV_X^D} \OO_X$ has automatically dimension $n=\dim X$.

Any Koszul free divisor $D$ is Spencer since, from the very definition, the complex $\DD_X \otimes_{\VV_X^D}\Sp^\bullet(\log D)$ is concentrated in cohomological degree $0$ (its graded complex with respect to $G^\bullet$ is concentrated in cohomological degree $0$), and $\DD_X\otimes_{\VV_X^D} \cO_X$ is holonomic since it is locally presented as $\DD_{X,p}/\DD_{X,p}\langle \delta_1,\dots,\delta_n\rangle$, where $\{\delta_1,\dots,\delta_n\}$ is a basis of $\fDer(-\log D)_p$, and the quotient of $\gr \DD_{X,p}$ by the ideal $\langle
\sigma(\delta_1),\dots,\sigma(\delta_n)\rangle$ has dimension $n$. Moreover, in this case $\gr \left( \DD_X\otimes_{\VV_X^D} \cO_X \right) \simeq  \gr \DD_X \otimes_{\gr \VV_X^D} \OO_X$, $\{\delta_1,\dots,\delta_n\}$ is an involutive basis of the ideal $\DD_{X,p}\langle \delta_1,\dots,\delta_n\rangle$ with respect to the order filtration and the characteristic variety of $\DD_X\otimes_{\VV_X^D} \cO_X$ is (locally) given by $\sigma(\delta_1)=\cdots=\sigma(\delta_n)=0$.
\medskip

Actually, the condition for a free divisor to be Koszul is equivalent to the fact that the logarithmic stratification of $D$ \cite[\S 3]{Saito-80} is locally finite, or equivalently, that any logarithmic stratum of $D$ is holonomic in the sense of {\em loc.~cit.}  (see \cite[Theorem (7.4)]{gmns}). In particular, any plane curve and any free hyperplane arrangement is a Koszul free divisor.

On the other hand, we know that any locally quasihomogeneous free divisor is Koszul \cite{cal_na_stek}\footnote{Actually, it is easy to see that for any locally quasihomogeneous divisor, free or not, the logarithmic stratification of $D$ is locally finite.}. Furthermore, the roots of the $b$-function of a reduced local equation $f=0$ of any locally quasihomogeneous free divisor are symmetric with respect to $-1$ (see \cite[Theorem 5.6]{cal_na_comp}, \cite[Corollary (4.2)]{nar_symmetry_BS}) and the $\DD[s]$-annhilator of $f^s$ is generated by order one operators \cite[Corollary 5.8]{cal_na_comp}. From there we deduce a purely algebraic proof of Theorem \ref{lct}, based on Theorem \ref{th:D-crit-LCT} (see \cite[Corollary (4.5) and Remark (4.6)]{nar_symmetry_BS}).

\begin{nota}
{Notice also that, for any hyperplane arrangement, free or not, with equation $f=0$, reduced or not, the $\cD_{X,p}$--annihilator of $f^{-1} \in \cO_{X,p}(*D)$ if generated by operators of order 1 \cite[Th. 5.3]{walther-2017}. As mentioned before, see Sect. \ref{conjecture-of-terao},  LCT  holds for any hyperplane arrangement \cite{Bath2022}.

If a free divisor $D\subset \CC^n$ is Spencer with a polynomial defining equation $f=0$, then by \cite[Crit. 3.1 and 4.1]{cas_ucha_exp}, $D$ satisfies LCT if and only if the annihilating ideal in $\cD$ of $1/f$ is generated by operators of order 1; that is, if condition (ii.2) in Remark \ref{lctconcrete} holds.
In \cite[Rk. 5.8]{cas_ucha_exp} examples are given of three free divisors in $\CC^3$, defined by quasi-homogeneous polynomials (with strictly positive weights),  that do not satisfy LCT and hence, by Th. \ref{lct}, they are not locally quasi-homogeneous. That gives a negative answer to a question proposed in \cite[Prob. 6.5]{cal_na_comp}.
}
\end{nota}

\begin{nota}
The sheaf $\check{\Omega}^\bu_D$ introduced in \cite{dfafd} (see \ref{sec:omega-cech}) is the quotient of $\Omega^\bullet_X$ by the subcomplex $\Omega^\bullet_X(\log D)(\cO_X(-D))$, and so there is a commutative diagram of complexes of sheaves of $\CC$-vector spaces
\beq \label{eq:big-cd}
\begin{tikzcd}
0 \ar[]{r}  & j_! \CC_{X\ssm D} \ar[]{r}{\text{adj.}} \ar[]{d}{\lambda} &  \CC_X \ar[]{r} \ar[]{d}{\text{q-i}}& \CC_D\ar[]{r} \ar[]{d}{\overline{\lambda}}& 0\\
0 \ar[]{r} & \Omega^\bullet_X(\log D)(\cO_X(-D)) \ar[]{r}{\text{incl.}} &  \Omega^\bullet_X \ar[]{r}& \check{\Omega}^\bu_D \ar[]{r}& 0
\end{tikzcd}
\eeq
with exact rows,
where the vertical arrows are the natural ones. The bottom row comes from the triangle
$$
\DD_X \Lotimes_{\VV_X^D} \OO_X  \longrightarrow \OO_X \longrightarrow K^\bullet \stackrel{+1}{\longrightarrow}
$$
by applying the $\DR(-)$ functor.

The middle vertical arrow in (\ref{eq:big-cd}) is a quasi-isomorphism by Poincaré lemma, and $\lambda$ is a quasi-isomorphism if and only if $\overline{\lambda}$ is a quasi-isomorphism. So, we can add another equivalent property to (2) in Theorem \ref{th:D-crit-LCT}:
\begin{enumerate}\em
\item[(2-4)] The canonical map $\CC_D \longrightarrow \check{\Omega}^\bu_D$ is a quasi-isomorphism.
\end{enumerate}

We deduce a new proof of Theorem \ref{lct} by using
\cite[Lemma 3.3]{dfafd}, where it is proven that the map $\overline{\lambda}$ is a quasi-isomorphism for any locally quasihomogeneous divisor (not necessarily free).
\end{nota}

In fact, a similar argument shows that LCT holds for any locally weakly quasihomogeneous free divisor and {\em a posteriori} any such free divisor is Spencer (see Definition \ref{def:LWQH})
(see \cite[Remark 3.11]{wqh} and \cite[Remark 1.7.4]{nar_comp_08}, although the statement in \cite[Remark 3.11]{wqh} that any logarithmic differential form has a non negative weight is wrong).
Namely, we have the following.

\begin{teorema} \label{th:LCT-LWQH}
If $D\subset X$ is a locally weakly quasihomogeneous free divisor, then $D$ satisfies the logarithmic comparison theorem.
\end{teorema}

\begin{prueba} By Theorem  \ref{th:D-crit-LCT} (2-3), we need to prove that the complex
$\Omega^\bullet_X(\log D)(\cO_X(-D))_p$ is exact for any $p\in D$. Since $D$ is locally weakly quasihomogeneous, there are local coordinates $x_1,\dots,x_r,x_{r+1},\dots,x_n$ centered at $p$, $r>0$, weights $(w_1,\dots,w_r,0,\dots,0)$, $w_i> 0$ for all $i=1,\dots,r$,   and a reduced local equation $f$ of $D$ at $p$ which is quasihomogeneous in these coordinates with $\wt(f)> 0$.

By the standard argument (see Lemma 3.3 in \cite{dfafd}), the complex
$$\Omega^\bullet_X(\log D)(\cO_X(-D))_p = f \Omega^\bullet_X(\log D)_p$$ is homotopic to its weight zero subcomplex, and the theorem is a consequence of the following lemma.
\end{prueba}

\begin{lema} Let $f$ be a quasihomogeneous polynomial in $\CC\{x_{r+1},\dots,x_n\}[x_1,\dots,x_r]$, of strictly positive weight with respect to weights $(w_1,\dots,w_r,0,\ld,0)$ with $w_1,\ld,w_r$ all strictly positive. Then, for any $p>0$ and any non-zero quasihomogeneous logarithmic $p$-form $\omega$, we have $\wt \omega
> -\wt f$.
\end{lema}

\begin{prueba} Since $f\omega$ must be holomorphic, we know that  $\wt \omega
\geq -\wt f $.
We may assume $\p f/\p x_1\neq 0$ (i.e. $f$ effectively depends on $x_1$).

Suppose that $\omega$ is a non-zero logarithmic $p$-form, and that $\wt \omega  = -\wt f $. Then, $\alpha = f \omega$ is a non-zero holomorphic $p$-form with $\wt \alpha =0$. That means that $\alpha$ only depends on the variables $x_{r+1},\dots,x_n$ of weight 0:
$$\alpha = \sum_{r+1\leq i_1<\ld<i_p\leq n} \alpha_{i_1,\ld,i_p}(x_{r+1},\dots,x_n)dx_{i_1}\wedge\ld\wedge dx_{i_p}.$$
Reordering the variables, we may assume $\alpha_{r+1,\ld,r+p}\neq 0$.
\medskip

Since $\omega$ is logarithmic, the $p+1$-form
$ df \wedge \omega$ is holomorphic (this follows from the fact that $f\omega$ and $fd\omega$ are holomorphic). The coefficient of $dx_1\wedge dx_{r+1}\wedge\ld\wedge dx_{r+p}$ in $df\wedge \omega$, which must be holomorphic, is
$$\frac{\alpha_{r+1,\ld,r+p}}{f} \frac{\p f}{\p x_1}.$$
So $f$ divides $\alpha_{r+1,\ld,r+p}\frac{\p f}{\p x_1}$ in the ring of power series
$\CC\{x_1,\dots,x_r,x_{r+1},\dots,x_n\}$. However, $f$ and $\p f/\p x_1$ are quasihomogeneous polynomials in $\CC\{x_{r+1},\dots,x_n\}[x_1,\dots,x_r]$. By equating homogeneous parts, we deduce that $f$ also divides $\alpha_{r+1,\ld,r+p}\p f/\p x_1$ in $\CC\{x_{r+1},\dots,x_n\}[x_1,\dots,x_r]$. Since $\wt(\p f/\p x_1)=\wt(f)-\wt(x_1)$, this is impossible.
\end{prueba}

A natural question is to characterize free divisors for which LCT holds. This has been done in \cite{cmnc} for plane curves, and more generally in \cite[Cor. 1.8]{torrelli_2004} and \cite[Th. (4.7)]{nar_symmetry_BS} for Koszul free divisors. For general free divisors the following conjecture remains open \cite[Conjecture 1.4]{cmnc}:

\begin{conjetura} Let $D\subset X$ a free divisor. If LCT holds for $D$, then $D$ is {\em strongly Euler homogeneous}\index{Strongly Euler homogeneous}, i.e. for each $p\in D$ there is a local reduced equation of $D$ at $p$ and a germ of vector field $\chi$ at $p$, singular at $p$, such that $\chi(f)=f$.
\end{conjetura}

In \cite[\S 3, 4]{castro-ucha-proc-2005} the authors provide an infinite family of free divisors for which LCT does not hold. Any divisor of this family is defined by a polynomial $$f_{p,q}=(x_1^p-x_2^q)\prod_{i=3}^n(x_1x_i+x_2)$$ for $n\geq 3$ and $3\leq p < q$.

\subsection{The logarithmic Spencer complex revisited}\index{Logarithmic Spencer complex}

Recall that, see (\ref{log-Spencer-complex}), for a free divisor $D$ in a complex manifold $X$, the logarithmic Spencer complex $\Sp^\bullet(\log D)$ is a locally free resolution of the left $\VV_X^D$-module $\OO_X$ \cite[Th. 3.1.2]{calde_ens}. We are going to
describe the complex
$\DD_X \otimes_{\VV_X^D} \Sp^\bullet(\log D)$ in terms of a (local) basis $\{\delta_1,\ldots,\delta_n\}$ of the free $\cO_X$--module $\fDer(-\log D)$.
We can write $$[\delta_i,\delta_j] = \sum_{k=1}^n \alpha_k^{ij}\delta_k$$ for $1\leq i<j\leq n$ where $\alpha_k^{ij} \in \cO_X$ for $k=1,\ldots,n$.
We also have $$\cD_X \otimes_{\VV_X^D} \Sp^p(\log D) = \cD_X \otimes_{\VV_X^D} {\VV_X^D} \otimes _{\cO_X} \bigwedge^p \fDer(-\log D) \simeq \cD_X \otimes _{\cO_X} \bigwedge^p \fDer(-\log D) $$ for $p=0,\ldots,n$, and we also denote $\epsilon_{-p}$ the differential of this complex of left $\cD_X$--modules.  We denote this complex by $\Sp^\bullet_{\cD_X}(\log D)$.

For $p\in \NN$ we denote  $\Lambda_p:= \{(i_1,\ldots, i_p)\in \NN^p\, \vert\, 1\leq i_1<i_2 < \ldots < i_p \leq n\}$. The free left $\cD_X$--module $\cD_X \otimes _{\cO_X} \bigwedge^p \fDer(-\log D)  \simeq \bigoplus_{\ubi \in \Lambda_p} \cD_X \,e_\ubi$ has rank $\binom{n}{p}$ with basis $\{e_\ubi \, \vert \ubi \in \Lambda_p \}$, for $p=0, \ldots,n$. An isomorphism
from $\cD_X \otimes _{\cO_X} \bigwedge^p \fDer(-\log D)$ onto $\bigoplus_{\ubi \in \Lambda_p} \cD_X \,e_\ubi$ maps $$1\otimes \delta_{i_1} \wedge \cdots \wedge \delta_{i_p} \,\mapsto\,  e_{\ubi}.$$

For $p\in \NN$ and $\ubi \in \Lambda_p$ we fix the following notations:
\begin{itemize}
\item The support of $\ubi$ is $\supp(\ubi) = \{i_1,\ldots,i_p\}$.
\item $\overline{\supp(\ubi)}=\{1,\ldots,n\}\,\, \ssm \,\, \supp(\ubi)$.
\item $\ubi(\widehat{k}) := (i_1,\ldots,i_{k-1}, i_{k+1}, \ldots,i_p)\in \Lambda_{p-1},$ for   $1\leq k \leq p$.
\item $\ubi(\widehat{k,\ell}) := \ubi(\widehat{k})(\widehat{\ell-1}) \in \Lambda_{p-2},$ $1\leq k <\ell  \leq p$.
\item $\sigma(q; \ubi) := \max\{j \in \{1,\ldots,p\}\, \vert \, i_j < q\}$ and $\sigma(q; \ubi)=0$ if $q < i_1$.
\item $\ubi(\widecheck{q}) = (i_1,\ldots,i_{k-1}, q, i_{k+1}, \ldots,i_p) \in \Lambda_{p+1}$ for $k= \sigma(q;\ubi)$ and $q\in\{1,\ldots,n\}\ssm \{i_1,\ldots,i_p\}$.
\end{itemize}

The complex $(\Sp^\bullet_{\cD_X}(\log D),\, \epsilon_{-\bullet})$ can be written as the complex $(\bigoplus_{\ubi \in \Lambda_\bullet} \cD_X \,e_\ubi,\, \widetilde{\epsilon_{-\bullet}})$ where
$$\widetilde{\epsilon_{-p}} : \bigoplus_{\ubi \in \Lambda_p} \cD_X \,e_\ubi \longrightarrow \bigoplus_{\ubj \in \Lambda_{p-1}} \cD_X \,e_\ubj$$ is the morphism of left $\cD_X$--modules defined by
$$\widetilde{\epsilon_{-p}} (e_\ubi) = \sum_{k=1}^p (-1)^{k-1} \delta_{i_k} e_{\ubi(\widehat{k})} + \sum_{1\leq \ell < m \leq p} (-1)^{\ell + m} \sum_{q \in \overline{\supp(\ubi(\widehat{\ell, m}))}} (-1)^{\sigma(q; \ubi(\widehat{\ell, m}))}\alpha^{i_\ell i_m}_q e_{\ubi(\widehat{\ell,m})(\widecheck{q})}.$$
Notice that the complex $(\bigoplus_{\ubi \in \Lambda_\bullet} \cD_X \,e_\ubi,\, \widetilde{\epsilon_{-\bullet}})$ can be viewed as a non commutative version of a Koszul complex.

Let us write down an example for $n=3$.

One has:

\begin{equation*}
\begin{split}
\epsilon_{-1}(e_i) & =\delta_i {\mbox{ for }} i=1,2,3, \\
\epsilon_{-2}(e_{ij}) & = \delta_i e_j - \delta_j e_i - \sum_{q=1}^3 \alpha_q^{ij} e_q \, {\mbox { for }}  1\leq i <j \leq 3, \\
\epsilon_{-3}(e_{123}) & = \delta_1 e_{23} - \delta_2 e_{13} + \delta_3 e_{12} + (\alpha_1^{13}+\alpha_2^{23}) e_{12} + (-\alpha_1^{12} + \alpha_3^{23})e_{13} + (-\alpha_2^{12} - \alpha_3^{13})e_{23} \\ & = (\delta_3+ \alpha_1^{13}+\alpha_2^{23}) e_{12} + (-\delta_2 -\alpha_1^{12} + \alpha_3^{23}) e_{13} + (\delta_1-\alpha_2^{12} - \alpha_3^{13}) e_{23}.
\end{split}
\end{equation*}

\subsection{The annihilator ideal of $1/f$ for a reduced equation $f$ of a normal crossing divisor in $\CC^n$}

Write   $\cD = \cD_{\CC^n,0}$ and $f = x_1\cdots x_r$ for $1\leq r \leq n$. One has the following equality  $$Ann_\cD(1/f) = \cD(\partial_{r+1}, \ldots, \partial_n, x_1\partial_1 +1, \ldots, x_r \partial_r +1).$$ Since the inclusion $A_n:=A_n(\CC) \subset \cD$ is flat, it is enough to prove a similar equality for the Weyl algebra $A_n$ instead of $\cD$.
The inclusion $$ A_n (\partial_{r+1}, \ldots, \partial_n, x_1\partial_1 +1, \ldots, x_r \partial_r +1) \subset \Ann_{A_n}\left(\frac{1}{f}\right)$$ is obvious since each generator of the first ideal annihilates $1/f$.
To prove the opposite inclusion, let us  assume first $1=r=n$ and write $x=x_1$ and $\partial = \partial_{x}$. First of all, any $P=P(x,\partial) \in A_1$ can be written as $$P=Q(x\partial + 1) + R + S$$ for unique $Q\in A_1$, $R:=R(x)\in \CC[x]$ and $S:=S(\partial)\in \CC[\partial](\partial)$. If $P(\frac{1}{x})=0$ and $S(\partial)$ is not zero with order $s\geq 1$, then $\frac{R}{x}+ S(\partial)(\frac{1}{x})=0$ which is a contradiction since the order of the pole of the rational function $\frac{R}{x}+ S(\partial)(\frac{1}{x})$ is $s+1$. Then $S=0$ and then $R=0$ and $P\in A_1(x\partial+1)$. The general case $1\leq r \leq n$ can be reduced to $r=n$ by taking, for any operator $P$ annihilating $\frac{1}{x_1\cdots x_r},$ the remainder of the division of $P$ by $\partial_{r+1},\ldots, \partial_n$. Finally, for the case $r=n$ we proceed  by division of an operator $P$ annihilating $\frac{1}{x_1\cdots x_n},$ with respect to $\partial_1 x_1, \ldots, \partial_n x_n$.

The $\cD$--module criterion for LCT gives another proof that LCT holds for normal crossing divisors. Indeed, any such divisor $D \equiv (x_1\cdots x_r=0) \subset \CC^n$ is Koszul free and then the complex
$\DD \otimes_{\VV}\Sp^\bullet(\log D)$ is concentrated in cohomological degree $0$ and $\DD \otimes_{\VV}\cO$ is locally presented as the quotient $\cD / \cD(\partial_{r+1}, \ldots, \partial_n, x_1\partial_1 +1, \ldots, x_r \partial_r +1)$, since  $\{x_1\partial_1, \ldots, x_r\partial_r,\partial_{r+1},\ldots, \partial_n\}$ is a basis of $\Der(-\log D)$. Finally, one has $$\frac{\cD}{ \cD(\partial_{r+1}, \ldots, \partial_n, x_1\partial_1 +1, \ldots, x_r \partial_r +1)} \simeq \cD f^{-1} = \cO(\star D).$$ The last equality follows from the well-know fact that the (global) Bernstein-Sato polynomial of $f=x_1\cdots x_r$ is $(s+1)^r$.

\noindent {\bf Acknowledgements}

The authors thank the editors of the Handbook of Geometry and Topology of
Singularities for the invitation to contribute this survey, and the referee and Alberto Castaño for very constructive comments.

FCJ acknowledges support from PID2020-117843GB-I00/AEI/10.13039/501100011033. LNM  acknowledges support from PID2020-114613GB-I00/AEI/10.13039/501100011033

\end{document}